\documentclass[11pt,a4paper,twoside]{article}

\usepackage[centertags]{amsmath}
\usepackage[T1]{fontenc}
\usepackage{amsfonts}
\usepackage{amssymb}
\usepackage{theorem}
\usepackage{ulem}
\usepackage{epsfig}
\usepackage{subfigure}
\usepackage[all]{xy}
\usepackage{changebar}
\usepackage{color}
\usepackage{mathrsfs}
\usepackage[boldsans]{ccfonts}
\usepackage{array}
\usepackage{tabulary}
\usepackage{supertabular}
\usepackage{calc}
\usepackage{eepic}
\usepackage{titlesec}
\usepackage{graphicx}
\usepackage{enumitem}

\theoremstyle{break} \newtheorem{theorem}{Theorem}[section]
\theoremstyle{break} \newtheorem{satz}{Theorem}
\theoremstyle{break} 
\theoremstyle{break} \newtheorem{definition}[theorem]{Definition} 
\theoremstyle{break} \newtheorem{lemma}[theorem]{Lemma}
\theoremstyle{break} \newtheorem{corollary}[theorem]{Corollary}
\theoremstyle{break} \newtheorem{example}[theorem]{Example}
\theoremstyle{break} 
{\theorembodyfont{\rmfamily}\newtheorem{remark}[theorem]{Remark}}
{\theorembodyfont{\rmfamily}}
\theoremstyle{break} 
\theoremstyle{break} 
\theoremstyle{break} 
\theoremstyle{break} 
\numberwithin{equation}{section}

\newcommand{\N}{{\mathbb{N}}}
\newcommand{\R}{{\mathbb{R}}}
\newcommand{\D}{{\mathbb{D}}}
\newcommand{\C}{{\mathbb{C}}}

\def\Re{\mathop{{\rm Re}}}

\newcommand{\hide}[1]{}

\begin{document}

\renewcommand{\thefootnote}{}
\stepcounter{footnote}


\smallskip

\begin{center}
{\Large \bf Critical sets of bounded analytic functions, zero sets of Bergman
spaces \\[2mm] and nonpositive curvature}
\footnote{2000 Mathematics Subject
 Classification: 30H05, 30J10, 35J60, 30H20, 30F45, 53A30\\ This research was
supported by the Deutsche Forschungsgemeinschaft (Grant: Ro 3462/3-1)}

\end{center}
\renewcommand{\thefootnote}{\arabic{footnote}}
\setcounter{footnote}{0}

\smallskip

\begin{center}
{\large Daniela Kraus}
\end{center}

\smallskip

\begin{center}
\begin{minipage}{13.5cm}
{\bf Abstract.}{ \small
A classical result due to Blaschke states that  for every analytic
self--map $f$ of the open unit disk of the
complex plane there exists a Blaschke product $B$ such that the zero sets of $f$ and $B$ agree. In this paper we show that there is an analogue statement for
critical sets, i.\!\;e.~for every analytic self--map $f$ of the open
unit disk there is even  an indestructible Blaschke product $B$ such that the critical
sets of $f$ and $B$ coincide. 
We further relate the problem of describing the critical sets of bounded
analytic functions to the problem of characterizing the zero sets of some
weighted Bergman space as well as to the Berger--Nirenberg problem from differential geometry. By
 solving the Berger--Nirenberg problem for a special case we identify the 
 critical sets of bounded analytic
 functions with the zero sets of the 
weighted Bergman space ${\cal A}_1^2$.
}
\end{minipage}
\end{center}

\smallskip
\section{Introduction}

A sequence $(z_j)$ of points  in a domain $G$ of the complex plane $\C$ is called the zero set of an analytic function
$f: G \to \C$, if $f$ vanishes precisely on this set $(z_j)$. This means that
$f(z)\not=0$ for $z\in G \backslash (z_j)$ and  if the
point $\xi \in G$
occurs $m$ times in the sequence $(z_j)$, then $f$ has a zero  at $\xi$ of precise
order $m$. 
By definition, the critical set of a nonconstant analytic function is the zero set of its first derivative. 
There is an extensive literature on critical sets.
In particular, there are many interesting results on the relation
between the zeros and the critical points of analytic and harmonic functions.
 A classical reference for this is the book of Walsh \cite{Wal1950}.

\smallskip

In this paper we study the problem of describing the critical sets of analytic
self--maps of the open unit disk $\D$ of the complex plane
$\C$. For this purpose
the classical characterization of the zero sets of bounded
analytic functions due to  Jensen \cite{Jen1899},
Blaschke \cite{Bla1915} and F.~and R.~Nevanlinna \cite{Nevs1922} serves
as a kind of model.

 \begin{satz}\label{thm:classical_zero}
 Let $ (z_j)$ be a sequence in  $\D$. Then the following statements are
equivalent. 
\begin{itemize}
 \item[(a)] There is an analytic self--map of $\D$ with zero set
 $ (z_j)$.
\item[(b)]  There is a Blaschke product $B$ with zero set  $
  (z_j)$, i.\;\!e.~$B(z)=\displaystyle \prod \limits_{j=1}^{\infty}\frac{\overline{z_j}}{|z_j|}\,
  \frac{z_j -z}{1-\overline{z_j}z}$.
\item[(c)] The sequence $(z_j)$ fulfills the Blaschke condition, i.\!\;e.~$
\sum \limits_{j=1}^{\infty} \big(1-|z_j|\big) < + \infty\, .
$

\item[(d)] There is a function in the Nevanlinna class  ${\cal N}$
 with zero set $(z_j)$.
\end{itemize}
\end{satz}

Let us recall that 
 a function $f$ analytic in $\D$  belongs to ${\cal
  N}$ if and only if the
integrals 
\begin{equation*}
\int \limits_{0}^{2\pi} \log^+|f(re^{it})|\, dt
\end{equation*}
 remain bounded as $r\to 1$.

\medskip

For the special case of a {\it finite} sequence a result related to Theorem
\ref{thm:classical_zero} but for critical
sets instead of zero sets can be found in  work of
Heins \cite[\S29]{Hei62}, 
Wang \& Peng \cite{WP79}, Zakeri \cite{Z96} and  Stephenson
\cite{Ste2005}: For
every finite sequence ${\cal C}=(z_j)$ in $\D$ there is always a
 {\it finite} Blaschke product whose critical set coincides with ${\cal C}$. 
A recent first generalization of this result to {\it infinite} sequences 
is discussed in \cite{KR}. There it is shown that every Blaschke sequence
$(z_j)$ is the critical set  of an infinite Blaschke product. However,
the converse to this, known as the Bloch--Nevanlinna conjecture \cite{Dur1969}, is false. There do exist
 Blaschke products, whose critical sets
fail to satisfy the Blaschke condition, see \cite[Theorem 3.6]{Col1985}. 
Thus the critical sets of bounded analytic functions are not
 just the Blaschke sequences and the situation for critical sets 
seems more subtle than for  zero sets.

\smallskip


The main result of this paper is
the following counterpart of
Theorem \ref{thm:classical_zero} for critical sets of bounded analytic functions.
\begin{theorem}\label{thm:main}
Let $ (z_j)$ be a sequence in $\D$. Then the following statements are
equivalent.
\begin{itemize}
 \item[(a)] There is an analytic  self--map of $\D$ with critical set
 $ (z_j)$.
\item[(b)]  There is an indestructible Blaschke product with critical set  $
  (z_j)$.
\item[(c)] There is  a function in the weighted Bergman space $ {\cal
    A}^2_1$ with zero set  $ (z_j)$.
\item[(d)] There is a function in the Nevanlinna class ${\cal N}$ with critical set $(z_j)$.
\end{itemize}
\end{theorem}

We note that a Blaschke product $B$ is said to be indestructible, if $T
\circ B$ is a Blaschke product
    for {\it every} unit disk automorphism $T$, see \cite[p.~51]{Col1985}.
The weighted Bergman space $ {\cal
    A}^2_1$ consists of all functions analytic in $\D$ for which  
\begin{equation*}
\iint \limits_{\D} (1-|z|^2) \, |f(z)|^{2}\,
d\sigma_z< + \infty\, ,
\end{equation*}
where $\sigma_z$ denotes two--dimensional Lebesgue measure with respect to
$z$, see for instance \cite[p.~2]{HKZ}.

\medskip

A few remarks are in order. 
First, implication (d) $\Rightarrow$ (a) of Theorem \ref{thm:main} is an old
result by Heins \cite[\S 30]{Hei62}. However, as part of this paper we
provide a new and different approach to this result.

\smallskip

Second, the simple geometric characterization of the zero sets of bounded analytic
functions via the Blaschke condition (c) in Theorem \ref{thm:classical_zero} has
not found an explicit counterpart for 
critical sets yet.  However, condition (c) of Theorem \ref{thm:main}  might be seen as
an implicit substitute.
The zero sets of (weighted) Bergman space functions have intensively been studied in the 1970's and
1990's by Horowitz
\cite{Hor1974,Hor1977}, Korenblum \cite{Kor1975} and Seip \cite{Seip1994,
  Seip1995}. As a result quite sharp necessary as well as sufficient
conditions for a sequence 
to be the zero set of a Bergman
space function are available.
In view of Theorem \ref{thm:main} all these results about zero sets of
Bergman space functions carry now
over to the critical sets of bounded analytic functions and vice versa.
Unfortunately, a {\it geometric}
characterization of the zero sets of  (weighted) Bergman space functions is
still unknown, ``and it is well known
that this problem is  very difficult'', cf.~\cite[p.~133]{HKZ}.

\medskip

Although, at first sight Theorem \ref{thm:main} appears
to be a result
exclusively in the realm of complex analysis, the ideas of its proof 
have their origin in differential geometry and partial differential equations.
We now give a brief account of the relevant interconnections.

\subsection*{Conformal metrics and associated analytic functions}

The proof of implication (a) $\Rightarrow$ (b) of Theorem \ref{thm:main} relies on conformal pseudometrics with Gauss curvature bounded
above by $-4$. The key steps are the following.
Suppose $f$ is a nonconstant analytic self--map of $\D$ with critical set
${\cal C}=(z_j)$. Then  
the pullback of the Poincar\'e metric
\begin{equation*}
\lambda_{\D}(w)\,|dw| = \frac{1}{1-|w|^2}\, |dw|
\end{equation*}
 for $\D$ with constant curvature $-4$ via $f$, i.\!\;e.~ 
\begin{equation*}
\lambda(z)\, |dz|=\frac{|f'(z)|}{1-|f(z)|^2}\,|dz|\, ,
\end{equation*}
induces a conformal pseudometric of constant curvature $-4$  which vanishes on
${\cal C}$ (cf.~Definition \ref{def:zeroset}). This allows us to apply a version of Perron's method\footnote{See Subsection
  \ref{sec:perron}, in particular Theorem \ref{thm:perron2}.} for conformal
pseudometrics which guarantees a  unique {\bf maximal} conformal pseudometric
$\lambda_{max}(z)\, |dz|$  on $\D$ with constant
curvature $-4$ which vanishes precisely on ${\cal C}$. 
Now an extension of Liouville's representation theorem (Theorem \ref{thm:liouville}) says
that every conformal pseudometric with constant curvature $-4$ and zero set
${\cal C}$  can be represented as the  pullback of the Poincar\'e metric
$\lambda_{\D}(w)\,|dw|$ via an analytic self--map of $\D$. In particular, 
\begin{equation*}
\lambda_{max}(z)\, |dz|= \frac{|F'(z)|}{1-|F(z)|^2}\,|dz|
\end{equation*}
for some analytic self--map $F$ of $\D$. 
We call $F$ a {\bf maximal function} 
for ${\cal C}$
since it is a ``developing map''
of the maximal conformal pseudometric with constant curvature $-4$
and zero set ${\cal C}$. Now, roughly speaking, the maximality of
$\lambda_{max}(z)\, |dz|$ forces its developing map $F$ to be maximal. More precisely, we have
  the following result.

\begin{theorem}\label{thm:0}
Every maximal function is an indestructible Blaschke product.
\end{theorem}

We note that in case of a finite sequence ${\cal C}$ the maximal functions are just the finite branched
coverings of $\D$.
One is therefore inclined to consider maximal functions for
infinite branch sets as ``infinite branched coverings'':

\begin{center}
\setlength{\extrarowheight}{6pt}
\begin{tabular}{|l|p{4.3cm}|l|}
\hline
 {\bf critical set} & { \bf maximal function} & {\bf mapping properties }\\[2mm] \hline
${\cal C}=\emptyset$ & automorphism of $\D$& unbranched covering of $\D$\!\;; \\[-2mm]
&&conformal self--map of $\D$\\[2mm] \hline
${\cal C}$ finite &finite Blaschke product &finite branched covering of $\D$\\[2mm] \hline
${\cal C}$ infinite & indestructible infinite \,\, Blaschke product& ``\!\;infinite
branched covering of $\D$\!\;''\\[7mm] \hline
\end{tabular}
\end{center}

\smallskip 

Hence maximal conformal pseudometrics of constant negative curvature and their associated maximal functions are
 of special interest for function theoretic con\-sider\-ations. 
The class of maximal conformal pseudometrics and their corresponding maximal
functions have already been studied by Heins in \cite[\S25 \& \S 26]{Hei62}. 
Heins discusses
  some necessary as well as sufficient conditions   
for maximal functions regarding their topological properties. He also posed the 
problem of characterizing maximal functions, cf.~\cite[\S 26, \S 29]{Hei62}. 
Theorem \ref{thm:0} gives  a partial answer to Heins' question.

\subsection*{The Gauss curvature PDE and the Berger--Nirenberg problem}

Heins' proof of implication (d)
$\Rightarrow$ (a) of Theorem \ref{thm:main}  relies
on conformal pseudometrics with curvature bounded above by $-4$.
Our approach is more general and will be based on the
Gauss curvature equation and an extension of Liouville's representation theorem, see Theorem
\ref{thm:liouville}. It uses the following idea.
Suppose $g$
is a nonconstant holomorphic function on $\D$. We consider the Gauss curvature equation
\begin{equation}\label{eq:gauss}
\Delta u = 4 \,|g'(z)|^2 e^{2u}
\end{equation}
on $\D$. Here $\Delta$ denotes the standard Laplacian. If we can guarantee the
existence of a realvalued $C^2$--solution
$u: \D \to \R$ to this PDE, 
then  $\lambda(z)\, |dz |:=|g'(z)|\,e^{u(z)}\, |dz|$ turns out to be  a  conformal pseudometric of constant
curvature $-4$ on $\D$ which vanishes on the critical set of $g$. Hence by Theorem \ref{thm:liouville}, 
\begin{equation*}
\lambda(z)\, |dz|=\frac{|f'(z)|}{1-|f(z)|^2}\, |dz|
\end{equation*}
for some analytic self--map $f$ of $\D$ and so the critical set of $f$ agrees
with the critical set of $g$. In short, if for
an analytic function $g$ in $\D$ the equation (\ref{eq:gauss}) has a
solution on $\D$, then there exists an analytic self--map $f$ of $\D$ such that
the critical sets of $g$ and $f$ coincide. Thus the main
point is  to characterize those holomorphic 
functions $g: \D \to \C$ for which  the PDE (\ref{eq:gauss}) has a solution on
$\D$.
In fact this problem is a special case of the well--known
Berger--Nirenberg problem from differential geometry, 
i.\!\;e.~the question,
whether for a Riemann surface $R$ and  a given  function $\kappa: R \to \R$
there exists a conformal metric on $R$ with Gauss curvature
$\kappa$.

\smallskip

We note that the Berger--Nirenberg problem is well--understood for the projective plane,
see \cite{Mos1973}, and has been extensively studied for compact Riemann surfaces,
see for instance \cite{Aub1998,Chang2004,Kaz1985,Str2005} as well as for the complex plane
\cite{Avi1986, CN1991, Ni1989}.  
However much less is known for proper subdomains $G$ of the complex plane,
see \cite{BK1986, HT92, KY1993}. In
this situation the Berger--Nirenberg problem
reduces  to the question if for a given function $k:G \to \R$ the Gauss curvature equation
\begin{equation}\label{eq:gauss1}
\Delta u= k(z)\, e^{2u}
\end{equation}
has a solution on $G$. We just note that $k$ is the
negative of the curvature $\kappa$ of the conformal metric $e^{u(z)}\, |dz|$.

In Theorem \ref{thm:sol1} we give some necessary as well as sufficient conditions
for the solvability of the Gauss curvature equation  (\ref{eq:gauss1})
only in terms of the  function $k$.
For instance, we shall see that the Gauss curvature equation (\ref{eq:gauss1}) has a solution  on $\D$ if $k$ is a nonnegative
locally H\"older continuous function on $\D$ and 
\begin{equation}\label{eq:sufcon}
\iint \limits_{\D} (1-|z|^2)\, k(z) \, d\sigma_z <+ \infty\, .
\end{equation}

\smallskip

 Let us make some remarks.
 First, this result generalizes previous results by Kalka $\&$ Yang in
 \cite{KY1993}. 
There the authors find a family of nonnegative locally H\"older continuous
functions $k_j$ on $\D$ tending uniformly to $+\infty$ at the boundary of $\D$
such that 
if $k$ is an essentially positive
function (see page \pageref{page:ess}) which satisfies $k\le k_j$ on $\D$ for some $j$,
then (\ref{eq:gauss1}) has a solution on $\D$.
 We note that all $k_j$ are radially symmetric and fulfill
(\ref{eq:sufcon}), see Example \ref{ex:k1} as well as the comments following
 Example \ref{ex:k1}.

\smallskip

Secondly, 
although condition (\ref{eq:sufcon}) is not necessary for the existence of a
solution to (\ref{eq:gauss1}) on $\D$ it is strong enough to deduce a necessary and
sufficient condition for the
solvability of the Gauss  curvature equation of the
particular form (\ref{eq:gauss}):

\begin{theorem}\label{thm:d}
Let $g:\D \to \C$ be an analytic function. Then the Gauss curvature equation
(\ref{eq:gauss}) has a solution on $\D$ if and only if $g'$ has
a representation as a product of an ${\cal A}_1^2$ function and a nonvanishing
analytic function in $\D$.
\end{theorem}

We note that Theorem \ref{thm:d} solves the Berger--Nirenberg problem for the
special case $R=\D$ and curvature functions of the form
$\kappa(z)=-|\varphi(z)|^2$ where $\varphi$ is analytic in $\D$. 

\medskip

A second observation is that   Theorem \ref{thm:d} 
leads  to a
characterization of the class of all holomorphic functions
$g: \D \to \C$ whose critical sets coincide with the critical sets of the class
of bounded analytic function, see Section \ref{sec:5} for more details.

\medskip

Finally, let us return to the implication (d) $\Rightarrow$ (a) of Theorem
\ref{thm:main}. Suppose $g$ is a function in  the Nevanlinna
class. Then  it turns out that $g'$ is the product of an ${\cal A}_1^2$ function and a
nonvanishing function, see the proof of Theorem \ref{thm:main} in Section \ref{sec:5}.
Consequently, by Theorem \ref{thm:d} the Gauss curvature equation
(\ref{eq:gauss}) has a solution which as we have seen implies that there
exists a bounded analytic function such that the critical
sets of these two functions agree.


\bigskip

We now will give a brief outline of the paper and  record in passing some
further results, which might be of interest in their own right. In Section \ref{sec:Metrics}  we discuss 
the theory of conformal pseudometrics  with curvature bounded above by $-4$ as
far as it is needed for this paper. 
We begin with some introductory material in Subsection \ref{sec:pseudometrics}.
Subsection \ref{sec:liouville} is devoted to Liouville's theorem and some of
its extensions.
We then relate in Subsection \ref{sec:boundary} the growth of a conformal pseudometric
with constant negative curvature  on $\D$  with inner functions.
This leads for instance to the following result, which might be viewed as an extension of Heins' characterization of finite Blaschke
products \cite{Hei86}\footnote{Cf. Theorem \ref{thm:boundary} in Subsection \ref{sec:boundary}.}.
 
\begin{corollary} \label{cor:5}
Let $f : \D \to \D$ be an analytic function. Then the following statements are
equivalent.
\begin{itemize}
\item[(a)] $\angle \lim \limits_{z \to \zeta} \left( 1-|z|^2 \right)
 \displaystyle \frac{|f'(z)|}{1-|f(z)|^2}=1$  for a.\!\;e.~$\zeta \in \partial
 \D$.
\item[(b)] $f$ is an inner function with finite angular derivative at almost
  every point of $\partial \D$.
\end{itemize}
\end{corollary}
Here, $\angle \lim $ denotes the nontangential limit.

\smallskip

Subsection \ref{sec:sk-metrics} will focus on some properties of
conformal pseudometrics with curvature bounded above by  $-4$.
Finally, in Subsection \ref{sec:perron}  we apply Perron's method to
guarantee the existence of maximal conformal pseudometrics with constant
curvature $-4$ and preassigned zeros.

\smallskip

In Section \ref{sec:Berger_Nirenberg_problem} we turn to the Berger--Nirenberg
problem for planar domains. Subsection \ref{sec:results} contains the main
results, illustrative  examples and remarks.
In particular, we establish  some necessary
and some sufficient conditions for the solvability of the Gauss curvature
equation only in terms of the curvature function and the domain in Theorem
\ref{thm:sol1} and Corollary \ref{cor:solutions1}. 
As a consequence of these conditions, we obtain
Theorem \ref{thm:d}.  The proofs of
the results are deferred to Subsection \ref{sec:proofs1} and Subsection \ref{sec:proofs2}.

\smallskip

Section \ref{sec:maximal} treats maximal functions. We begin with the proof of Theorem
\ref{thm:0} as well as of the equivalence of (a) and (b) in Theorem
\ref{thm:main}. We then discuss  maximal functions whose critical sets form
finite and Blaschke sequences respectively in more detail.
For instance,
maximal functions with finite critical
sets are finite Blaschke products and vice versa.
This together with Theorem \ref{thm:0}  implies a series of refined Schwarz--Pick type inequalities:

\begin{corollary}\label{cor:3}
Let $f: \D \to \D$ be a nonconstant analytic function with critical set ${\cal
  C}$ and let ${\cal C}^*$ be a subsequence of ${\cal C}$. 
Then there exists a (indestructible) Blaschke
product $F$ with critical set ${\cal C}^*$
such that
\begin{equation*}
\frac{|f'(z)|}{1-|f(z)|^2} \le \frac{|F'(z)|}{1-|F(z)|^2}\, ,\quad  z \in \D\!\;.
\end{equation*}

If ${\cal C}^*$ is finite, then $F$ is a finite Blaschke product.

\smallskip

Furthermore, $f=T \circ F$ for some automorphism $T$ of $\D$
if and only if
\begin{equation*}
\lim_{z \to w} \frac{|f'(z)|}{1-|f(z)|^2} \,
\frac{1-|F(z)|^2}{|F'(z)|}=1
\end{equation*}
for some $w \in \D$.
\end{corollary}

We note that   for {\it finite}  sequences ${\cal    C}^*$ Corollary
\ref{cor:3} is a classical result of Nehari \cite{Neh1946}. 
Since {\it infinite} sequences are explicitly allowed in Corollary \ref{cor:3}, 
it generalizes Nehari's result.
In addition, for ${\cal C}^*={\cal C} \not= \emptyset$ we get a 
best possible sharpening
\begin{equation*}
\frac{|f'(z)|}{1-|f(z)|^2} \le \frac{|F'(z)|}{1-|F(z)|^2} < \frac{1}{1-|z|^2}
\end{equation*}
of the Schwarz--Pick inequality
\begin{equation*}
\frac{|f'(z)|}{1-|f(z)|^2} \le \frac{1}{1-|z|^2} \, .
\end{equation*}
We end Section \ref{sec:maximal} with a criterion for maximal functions with
finite and Blaschke sequences as critical sets.

\smallskip

In a short final Section \ref{sec:5} we prove  Theorem \ref{thm:main} and conclude
by some additional remarks.

\bigskip

Before beginning with the details it is worth making some comment about
notation. The action in this paper takes place on domains in the complex
plane. The letters $D$ and $G$ exclusively denote planar domains and will be
used without further explanation.


\section{Conformal  metrics and pseudometrics}
\label{sec:Metrics}
The following subsections discuss some selected topics on conformal
pseudometrics with negative curvature.  For more information we refer to
\cite{BM2007, Hei62, KL2007, KR2008, Smi1986}.

\subsection{Conformal metrics and pseudometrics}
\label{sec:pseudometrics}

We begin our brief account of conformal metrics and pseudometrics with some
basic definitions and results.

\medskip

\begin{definition}
A nonnegative continuous function $\lambda$ on
$G$, $\lambda: G \to [0, + \infty)$, $\lambda \not\equiv 0$, is
called conformal density on $G$ and the corresponding
 quantity \label{def:leng}
$\lambda(z) \, |dz|$   conformal pseudometric on $G$.
If $\lambda(z)>0$ for all $z \in G$, we say $\lambda(z) \, |dz|$ is a
conformal metric on $G$. We call a conformal
pseudometric $\lambda(z) \, |dz|$  regular on $G$, if $\lambda$ is 
of class $C^2$
 in $\{z \in G \, : \, \lambda(z)>0\}$\footnote{$C^2(G)$ denotes the set of realvalued twice
  continuously differentiable functions on $G$.}.
\end{definition}

We wish to emphasize that, according to our definition, $\lambda\equiv 0$ is {\bf not} a
conformal density (of a conformal pseudometric).

\smallskip 

A second remark is that some authors call a nonnegative upper semicontinuous
function a conformal density. For our applications however it suffices to ask
for continuity.

\smallskip

A geometric quantity associated with a conformal pseudometric
is its Gauss curvature.

\begin{definition}[Gauss curvature]\label{def:curvature}
Let $\lambda (z) \, |dz|$ be a regular conformal pseudometric on $G$. Then 
the (Gauss) curvature $\kappa_{\lambda}$ of $\lambda(z) \, |dz|$ 
is defined by
\begin{equation*}
\kappa_{\lambda}(z):=-\frac{\Delta (\log \lambda)(z)}{\lambda(z)^2}
\end{equation*}
 for all points $z \in G$ where $\lambda(z) > 0$. 
\end{definition}

An important property of the Gauss curvature is its conformal invariance. It
is based on the following definition.

\begin{definition}[Pullback of conformal pseudometrics]
Let $\lambda(w) \, |dw|$ be a conformal pseudometric
on $D$  and $w=f(z)$ be  a nonconstant analytic map from $G$ to $D$. Then
the conformal pseudometric
\begin{equation*} 
 (f^*\lambda)(z) \, |dz|:=\lambda(f(z)) \, |f'(z)| \, |dz|
\end{equation*}
defined on $G$, is called the pullback of $\lambda(w) \, |dw|$ under the map $f$.
\end{definition}

\begin{theorem}[Theorema Egregium]
For every analytic map $w=f(z)$ and every regular conformal pseudometric $\lambda(w)\,|dw|$ the relation
\begin{equation*} 
\kappa_{f^*\lambda}(z)=\kappa_{\lambda}(f(z))
\end{equation*}
is satisfied provided $\lambda(f(z)) \,  |f'(z)|>0$. 
\end{theorem}

Definition \ref{def:curvature} shows that if $\lambda(z)\, |dz|$ is a regular conformal metric
with curvature $\kappa_{\lambda}=\kappa$ on $G$, then the function $u:= \log \lambda$ 
satisfies the partial differential equation
\begin{equation}\label{eq:pde1}
\Delta u=-\kappa(z) \, e^{2\;\!u}
\end{equation}
on $G$. If, conversely, a $C^2$--function $u$ fulfills (\ref{eq:pde1}) on $G$,
then $\lambda(z):=e^{u(z)}$ induces a regular conformal metric on $G$ with
curvature $\kappa_{\lambda}=\kappa$.

\medskip

The ubiquitous example of a conformal metric is the Poincar\'{e} or hyperbolic
metric $\lambda_{\D}(z)\, |dz|$ for the unit disk $\D$  with constant
curvature $-4$. It has the following important property.

\begin{theorem}[Fundamental Theorem]\label{thm:ahlfors}
Let $\lambda(z)\,|dz|$ be a regular conformal pseudometric on
$\D$ with curvature bounded above by $-4$. Then $\lambda(z) \le
\lambda_{\D}(z)$ for every $z \in \D$.
\end{theorem}

Theorem \ref{thm:ahlfors} is due to Ahlfors \cite{Ahl1938} and it is usually  
called Ahlfors' lemma. However, in view of its relevance Beardon and Minda
proposed to call Ahlfors' lemma the  {\bf fundamental theorem}. We
will follow  their suggestion in this paper.

\subsection{Liouville's Theorem}
\label{sec:liouville}

Conformal pseudometrics of constant curvature $-4$ have a special nature.
First, they give us via the pullback a means of constructing conformal
pseudometrics with various prescribed properties without changing their
curvature. Second, every conformal pseudometric of constant curvature $-4$
can  locally be
represented by a holomorphic function. This is Liouville's theorem. In order
to give a precise statement we begin with a formal definition. 

 \begin{definition}[Zero set]\label{def:zeroset}
Let $\lambda(z)\, |dz|$ be a conformal pseudometric on $G$. We say
$\lambda(z)\, |dz|$ has a zero of order $m_0>0$ at $z_0 \in G$ if
\begin{equation*}
\lim_{z \to z_0} \frac{\lambda(z)}{|z-z_0|^{m_0}} \quad\text{ exists and }
\not=0 \, . 
\end{equation*}

\smallskip

We will call a  sequence ${\cal C}=(\xi_j) \subset G$
\begin{equation*}
(\xi_j):=(\underbrace{z_1, \ldots, z_1}_{m_1  -\text{times}},\underbrace{z_2,
  \ldots, z_2}_{m_2-\text{times}} , \ldots  )\,, \, \,  z_k \not=z_n \text{ if
} k\not=n, 
\end{equation*}
  the zero set of a conformal
pseudometric $\lambda(z) \, |dz|$, if $\lambda(z)>0$ for $z \in
G\backslash {\cal C}$ and if $\lambda(z)\, |dz|$ has a zero of order $m_k \in \N$
at $z_k$ for all $k$.

\end{definition}

\begin{theorem}[Liouville's Theorem]
\label{thm:liouville}

Let ${\cal C}$ be a sequence of points in a simply connected domain $G$
and let
$\lambda(z) \, |dz|$ be a regular conformal pseudometric on $G$ with constant curvature
$-4$ on $G$ and zero set ${\cal C}$. Then $\lambda(z) \, |dz|$ is the pullback of the
hyperbolic metric $\lambda_{\D}(w)\, |dw|$ under some analytic map $f:G \to
\D$, i.\!\;e.
\begin{equation}\label{eq:liouville}
 \lambda(z) =\frac{|f'(z)|}{1-|f(z)|^2}\, , \quad z \in G. 
\end{equation}
If $g:G  \to \mathbb{D}$ is another analytic function, then
$\lambda(z)=(g^*\lambda_{\D})(z)$ for all $z\in G$
 if and only if $g=T\circ f$ for some automorphism $T$ of $\D$. 
\end{theorem}

A holomorphic function $f$ with  property (\ref{eq:liouville})
 will be called a {\bf
  developing map} for $\lambda(z) \, |dz|$.

\medskip

Note that the critical set of each developing map coincides with the zero set
of the corresponding conformal pseudometric.  

\medskip

For later applications we wish to mention the following variant of Theorem \ref{thm:liouville}.

\begin{remark}\label{rem:liouville}
Let $G$ be a simply connected domain and let $\varphi:G \to \C$, $\varphi \not
\equiv 0$, be an analytic
map. If $\lambda(z)\, |dz|$ is a regular conformal metric with curvature $-4
\, |\varphi(z)|^2$, then there exists a holomorphic function $f:G \to \D$ such
that
\begin{equation*}
 \lambda(z) =\frac{1}{|\varphi(z)|}\frac{|f'(z)|}{1-|f(z)|^2}\, , \quad z
\in G. 
\end{equation*}
Moreover, $f$ is uniquely determined up to postcomposition with a unit disk automorphism.
\end{remark}

 Liouville \cite{Lio1853} stated
Theorem \ref{thm:liouville}   for the special case that $\lambda(z)
\, |dz|$ is a regular conformal metric.
We therefore like to refer to
 Theorem \ref{thm:liouville} as well as to Remark \ref{rem:liouville} as
 Liouville's theorem. 

\smallskip

 Theorem \ref{thm:liouville} and in particular the special case that $\lambda(z) \, |dz|$ is
 a conformal metric has a number of different proofs, see for instance 
\cite{Bie16, CW94, CW95, 
Min, Nit57, Yam1988}. Remark
\ref{rem:liouville} is discussed in \cite{KR}.

\subsection{Boundary behavior of developing maps}
\label{sec:boundary}

By Liouville's theorem it is perhaps not too surprising that there is some
relation between the  boundary
behavior of a conformal pseudometric  and the boundary behavior of a corresponding
developing map. The next result illustrates this relation.

\smallskip

\label{page1}
\begin{satz}[\text{\small cf.~\cite{Hei86, KRR06}\!\;}] \label{thm:boundary}
Let $f: \D \to \D$ be an analytic function and $I \subset \partial \D$ some
open arc. Then the following are equivalent.
\begin{itemize}
\item[(a)] 
$ \lim \limits_{z \to \zeta} \left( 1-|z|^2 \right)
 \displaystyle \frac{|f'(z)|}{1-|f(z)|^2}=1 \qquad \text{ for every } 
\zeta \in I \, , $
\item[(b)]$ \lim \limits_{z \to \zeta} 
 \displaystyle \frac{|f'(z)|}{1-|f(z)|^2}=+\infty \qquad \text{ for every } 
\zeta \in I \, , $
\item[(c)]
$f$ has a holomorphic extension across the arc $I$ with $f(I) \subset \partial \D$.
\end{itemize}
In particular, if $I=\partial \D$, then $f$ is a finite Blaschke product.
\end{satz}

\medskip

In fact similar results can be derived when the unrestricted
limits are replaced by angular limits.

\begin{lemma} \label{lem:lemma2}
Let $f : \D \to \D$ be an analytic function and $I$ some subset of $\partial\D$.
\begin{itemize}
\item[(1)]
If 
\begin{equation*}
\angle \lim \limits_{z \to \zeta} \left( 1-|z|^2 \right)
 \displaystyle \frac{|f'(z)|}{1-|f(z)|^2}=1 \qquad \text{ for every } 
\zeta \in I \, , 
\end{equation*}
then  $f$ has a finite angular derivative\footnote{see \cite[p.~57]{Sha1993}.} at a.\!\;e.~$\zeta
  \in I$. In particular,  
\begin{equation*}
 \angle \lim \limits_{z \to \zeta} |f(z)|=1 \qquad \text{ for a.\!\;e. } \zeta \in  I
\, .
\end{equation*}

\item[(2)]
If $f$ has a finite angular derivative (and $\angle \lim_{z \to \zeta} |f(z)|=1$) at some $\zeta
  \in I$, then
\begin{equation*}
\angle \lim \limits_{z \to \zeta} \left( 1-|z|^2 \right)
 \displaystyle \frac{|f'(z)|}{1-|f(z)|^2}=1 \, .
\end{equation*}
\end{itemize}
\end{lemma}

\medskip

Now, if $I=\partial \D$, then as a direct consequence of Lemma \ref{lem:lemma2} we obtain 
Corollary \ref{cor:5} in the Introduction.

\medskip

{\bf Proof of Lemma \ref{lem:lemma2}.} \\
(1) Let $\zeta \in I$
and assume that $\liminf \limits_{z \to \zeta} \displaystyle
\frac{1-|f(z)|}{1-|z|}=+\infty$. Since
\begin{equation*} |f'(z)|=\left( 1-|z|^2 \right) \frac{|f'(z)|}{1-|f(z)|^2}
\frac{1-|f(z)|^2}{1-|z|^2} \, , 
\end{equation*}
we deduce 
\begin{equation*}
 \angle \lim \limits_{z \to \zeta} |f'(z)|=+\infty \, . 
\end{equation*}

By Privalov's theorem, cf.~\cite[p.~47]{Pri1956}, this is only possible for $\zeta$ from a nullset
$I' \subset \partial \D$.

Thus
\begin{equation*}
\liminf \limits_{z \to \zeta} \displaystyle
\frac{1-|f(z)|}{1-|z|}<+\infty \qquad \text{ for a.\!\;e. } \zeta \in  I
\, .
\end{equation*}
Therefore,  $f$ has a finite angular derivative $f'(\zeta)$ at a.\!\;e.~$\zeta
\in I$, see \cite[p.~57]{Sha1993}.

\medskip

(2) For convenience we may assume $\zeta=1$ and $f(1):=\angle \lim_{z \to 1}
  f(z)=1$. Thus $f'(1):=\angle \lim_{z \to 1} f'(z)>0$, see
  \cite[p.~57]{Sha1993}. We define for $z \in \D$
\begin{equation*}
\varrho(z):= \frac{f(z)-1}{z-1} -f'(1)\, .
\end{equation*}
Then $\angle \lim_{z \to 1}\varrho(z)=0$ and $f(z)=1+f'(1)\, (z-1) +
\varrho(z)\, (z-1)$ for $z \in \D$.
Hence we can write
\begin{equation*}
1-|f(z)|^2= 2\, f'(1)\, \Re(1-z)+ \Re(1-z)\,r(z)\, , 
\end{equation*}
where $\angle \lim_{z \to 1}r(z)=0$. This yields
\begin{equation*}
 \left( 1-|z|^2 \right)
  \frac{|f'(z)|}{1-|f(z)|^2}= \frac{1-|z|^2}{2\;\!\Re(1-z)}\,
    \frac{|f'(z)|}{f'(1)}\, \frac{1}{1+\frac{r(z)}{2\, f'(1)}}\, .
\end{equation*}
If we now choose $z \in S_{\delta}:=\{\, 1+r\, e^{i\, \alpha} \in \D \, : \,
r>0, \, 
\, 
\alpha\in [ \pi/2+ \delta, 3\pi/2 - \delta]\, \} $ for $\delta >0$, then
we have $1-|z|^2=-2\;\! r \cos \alpha -r^2$ and $\Re(1-z)=-r\, \cos \alpha$.
Hence
\begin{equation*}
 \angle \lim \limits_{z \to 1}\frac{1-|z|^2}{2\;\!\Re(1-z)}=1+  \angle \lim \limits_{z \to 1} \frac{r}{2\cos\alpha}\ge 1+  \lim \limits_{r \to 0} \frac{r}{2\,
  \cos(\frac{\pi}{2} +\delta)}=1\, .
\end{equation*}
So we can conclude that
\begin{equation*}
 \angle \lim \limits_{z \to \zeta} \left( 1-|z|^2 \right)
  \frac{|f'(z)|}{1-|f(z)|^2}\ge 1\,.
\end{equation*}
Combining this with the Schwarz--Pick lemma gives the desired result.
\hfill{$\blacksquare$}

\medskip

Finally, there is a counterpart of implication (b) $\Rightarrow$ (c) of 
Theorem \ref{thm:boundary}.

\begin{lemma} \label{lem:lemma}
Let $f : \D \to \D$ be an analytic function.
If 
\begin{equation}\label{eq:innerfunction}
 \angle \lim \limits_{z \to \zeta} \frac{|f'(z)|}{1-|f(z)|^2}=+\infty \,
\qquad \text{ for a.\!\;e. } \zeta \in \partial \D \, , 
\end{equation}
then $f$  is an inner function.
\end{lemma}

{\bf Proof.}
 Assume $f$ is not inner. Then the angular limit
$f(\zeta)$ of $f$ exists and belongs to $\D$ for every $\zeta \in I$ for a set
$I \subseteq \partial \D$ of positive measure. Now, in view of (\ref{eq:innerfunction}) the angular limit of $f'$ would
be $\infty$ for a set 
$I'\subseteq \partial \D$ of positive measure, contradicting Privalov's theorem.
\hfill{$\blacksquare$}

\bigskip

It would be interesting to see an example of an inner function for which 
(\ref{eq:innerfunction}) is not true. We further note that
conditions (1) and (2) in Lemma \ref{lem:lemma2} do not complement each
other. Therefore we may ask if an analytic self--map $f$ of $\D$ which 
satisfies
\begin{equation*}
\angle\lim \limits_{z \to 1} \frac{|f'(z)|}{1-|f(z)|^2}\, (1-|z|^2)=1  
\end{equation*}
does have
an angular limit or even a finite angular derivative  at
$z=1$; this might then be viewed as a  converse of the
Julia--Wolff--Carath\'eodory theorem, see \cite[p.~57]{Sha1993}.

\subsection{SK--metrics}
\label{sec:sk-metrics}

In the next two subsections we take a brief look at Heins' theory of 
SK--metrics\footnote{In particular \S 2, \S 3, \S 7, \S 10, \S 12 and  \S 13
  in \cite{Hei62}; see also \cite{KR2008}.}, refine and extend it
 in order to give a self--contained overview of the results which are needed
 for this paper. 
An SK--metric in the sense of Heins 
is a conformal pseudometric 
whose ``generalized curvature'' is bounded
above by $-4$. More precisely, this ``generalized curvature'' is obtained by
replacing the standard Laplacian in Definition \ref{def:curvature} by the  generalized
lower Laplace operator, which is defined for
a continuous function $u$ by
\begin{equation*}
\Delta u(z)= \liminf_{ r \to 0} \, \frac{4}{r^2} \left(
  \frac{1}{2 \pi} \int \limits_{0}^{2 \pi} u(z +r e^{it}) \, dt - u(z)    
 \right) \, .
\end{equation*}
We note that in case $u$ is a $C^2$--function  the generalized lower Laplace operator coincides with the standard Laplace operator.
Hence one can  assign to an arbitrary 
conformal pseudometric
$\lambda(z)\, |dz|$  a Gauss curvature $\kappa_{\lambda}$ in a natural way.

\begin{definition}[SK--metric]\label{def:sk_metric}
A conformal pseudometric $\lambda(z)\, |dz|$ on $G$  is called
{SK--metric} on $G$ if  $\kappa_{\lambda}(z)\le -4$ for all $z \in G$ where $\lambda(z)>0$.
\end{definition}

Note that every SK--metric is a subharmonic function. Second, if the curvature
of an SK--metric $\lambda(z)\, |dz|$ is locally H\"older continuous on $G$, then $\lambda$ is
regular on $G$ by elliptic regularity (cf.~Theorem \ref{thm:existence2}).

\smallskip

We now record some basic but essential properties of SK--metrics.

\begin{lemma} [\text{\small cf.~\cite[\S 10]{Hei62}}\!\;] \label{lem:maxmetrics}
Let $\lambda(z) \, |dz|$ and $\mu(z) \, |dz|$ be SK--metrics on $G$.
Then $\sigma(z) :=  \max \{\lambda(z), \mu(z)\}$ induces an SK--metric on $G$.
\end{lemma}

\begin{lemma}[\text{\small cf.~\cite[Lemma 3.7]{KR2008}}\!\;) \, (Gluing Lemma] \label{lem:gluing}
Let  $\lambda(z)\, |dz|$ be an SK--metric on $G$ and let $\mu(z)\, |dz|$
be an SK--metric on a subdomain  $D$ of $G$ such that
the ``gluing condition''
\begin{equation*} 
\limsup \limits_{ D \ni z \to \zeta} \mu(z) \le  \lambda(\zeta)  
\end{equation*}
holds for all $\zeta \in \partial D \cap G$. Then $\sigma(z)\, |dz|$ defined by
\begin{equation*}
\sigma(z):=\begin{cases} \,   \max \{\lambda(z), \mu(z)\}\, ,   & \hspace{2mm} \,
  \text{if } \,  z \in  D, \\[2mm]
                       \,        \lambda(z) \,  ,       & \hspace{2mm} \,
                       \text{if } \, z \in G \backslash D,
          \end{cases} 
\end{equation*}
is an SK--metric on $G$.
\end{lemma}

\begin{lemma}[\text{\small cf.~\cite[\S 10]{Hei62}}\!\;]\label{lem:new_sk}
Let $\lambda(z)\, |dz|$ be an SK--metric on $G$ and let $s$  be a
nonpositive  subharmonic function on $G$, then
$\mu(z)\, |dz|:=e^{s(z)}\,  \lambda(z)\, |dz|$ is an SK--metric on $G$.
\end{lemma}

\medskip

\begin{theorem}[\text{\small cf.~\cite[\S 2]{Hei62}}\!\;) \, (Generalized Maximum Principle]\label{thm:gmp}
Let $\lambda(z)\, |dz|$ be an SK--metric on $G$ and $\mu(z)\, |dz|$ be a
regular conformal
metric on $G$ with constant curvature $-4$. If
\begin{equation*}
\limsup\limits_{z \to \zeta} \frac{\lambda(z)}{\mu(z)} \le 1 \quad \text{for
  all } \zeta \in \partial_{\infty} G\, \footnote{$\partial_{\infty} G$ means the
  boundary of $G$ in $\C \cup \{ \infty\}$.}\, ,
\end{equation*}
then $\lambda(z) \le \mu(z)$ for $z \in G$.
\end{theorem}

The following lemma provides a converse to the generalized maximum principle.
It also might be viewed as an alternative definition of
an SK--metric.

\begin{lemma}[\text{\small cf.~\cite[\S 3]{Hei62}}\!\;]\label{lem:char_metrics}
Let $\lambda$ be a continuous function on  $G$. Then the following are equivalent.
\begin{itemize}
\item[(a)]
$\lambda(z) \, |dz|$ is an SK--metric on $G$.
\item[(b)]
Whenever $D$ is a relatively compact subdomain of $G$, and $\mu(z) \,
|dz|$ is a regular conformal metric with constant curvature $ -4$
on $D$ satisfying
\begin{equation*}
\limsup_{z \to \zeta}  \frac{\lambda(z)}{\mu(z)} \le 1
\end{equation*}
for all $\zeta \in \partial D$, then $\lambda(z) \le \mu(z)$ for $z \in D$.
\end{itemize}
\end{lemma}

We now apply Lemma \ref{lem:char_metrics} to prove a removable singularity
theorem for SK--metrics.

\begin{lemma}[Removable singularities]\label{lem:hebsing}
Let $\lambda$ be a continuous function on $G$ which induces an SK--metric on
$G\backslash \{ z_0\}$. Then $\lambda(z) \, |dz|$ is an SK--metric on $G$.
\end{lemma}

{\bf Proof.}
Let $D$   be a  relatively compact subdomain of $ G$ which contains $z_0$ and 
let $\mu(z)\, |dz |$ be a regular conformal metric with constant curvature
$-4$ on $D$ such that 
\begin{equation*}
 \limsup_{z \to \zeta} \frac{\lambda(z)}{\mu(z)} \le 1   \quad \text{for
  all } \zeta \in \partial D\, .
\end{equation*}
Then the nonnegative function
\begin{equation*}
s(z):= \log^+\left( \frac{\lambda(z)}{\mu(z)} \right)=
\max\left\{ 0,\, \log\left( \frac{\lambda(z)}{\mu(z)} \right) \right\}\, ,
\quad z \in D\,,  
\end{equation*}
is subharmonic on $D \backslash \{z_0 \}$.
To see this, let  $s(z_*)>0$ at some point $z_*\in D \backslash \{z_0 \} $. Thus
$s(z)>0$ in a neighborhood of $z_*$ and consequently $\Delta s (z) \ge
4\, (\lambda(z)^2-\mu(z)^2) >0$ in this neighborhood, i.\!\;e.~$s$ is subharmonic
there. If
$s(z_*)=0$ for some $z_* \in D \backslash \{z_0 \}$ then $s$ satisfies the
submean inequality 
\begin{equation*}
s(z_*)=0\le \frac{1}{2\,\pi} \int \limits_0^{2\pi} s(z_*+re^{it})\, dt 
\end{equation*}
for all small $r$. We note that $s$  has a subharmonic extension to $D$,
since $s$ is bounded near $z_0$. 
By hypothesis
$\limsup_{z \to \zeta } s(z)=0$ for all $\zeta \in \partial D$ and so the maximum
principle for subharmonic functions implies that $s\equiv 0$, i.\!\;e.~$\lambda
\le \mu$ in $D$. 
Finally, by Lemma
\ref{lem:char_metrics},  $\lambda(z) \, |dz|$ is an SK--metric
on $G$.~\hfill{$\blacksquare$}

\subsection{Perron's method}
\label{sec:perron}

Perron's method for subharmonic functions \cite{Per1923} is used to treat the
classical Dirichlet problem in arbitrary bounded domains. One attractive
feature of this method is that it separates the interior existence problem
from that of the boundary behavior of the solution. In addition the solution
is characterized by a maximality property. Perron's method, for instance, can
be imitated to ensure the existence of solutions to fairly general elliptic
PDEs, cf.~\cite[Chapter 6.3]{GT}. 
We apply Perron's method to guarantee the existence of maximal conformal SK--metrics with prescribed
zeros. In passing  let us quickly recall the definition of a  Perron family.
For more information regarding Perron families for SK--metrics, see
\cite[\S 12 \& \S 13]{Hei62} and \cite[Section 3.2]{KR2008}.

\begin{lemma}[Modification]\label{lem:mod}
Let $\lambda(z)\, |dz|$ be an SK--metric on $G$  and let $K$ be  an open disk which
is compactly contained in $G$. Then there exists a unique SK--metric
$M_K\lambda(z)\, |dz|$ on $G$, called modification of $\lambda$ on $K$, with the following properties:
\begin{itemize}
\item[(i)] $M_K\lambda(z)  =\lambda(z)$  for every $z \in G\backslash K$ and
$M_K\lambda(z) \ge \lambda(z)$ for every $z \in K$,
\item[(ii)] $M_K\lambda(z)\, |dz| $ is a regular conformal metric on $K$ with
 constant curvature $-4$. 
\end{itemize}
\end{lemma}

\begin{definition}[Perron family]
A  family $\Phi$ of (densities of) SK--metrics on $G$  is called a Perron
family, if the following conditions are satisfied: 
\begin{itemize}
\item[(i)]  If $\lambda \in \Phi$ and $\mu \in \Phi$, then $\sigma \in \Phi$,
  where $\sigma(z) := \max \{\lambda(z), \mu(z)\}$ for $z \in G$. 
\item[(ii)] If $\lambda  \in \Phi$, then $M_K\lambda \in
  \Phi$ for any open disk $K$ compactly contained in~$G$.
\end{itemize}
\end{definition}

\begin{theorem} \label{thm:perron0}
Let $\Phi$ be a Perron family of SK--metrics on $G$.
If $\Phi \not= \emptyset$, then 
\begin{equation*} 
\lambda_{\Phi}(z):=\sup_{\lambda \in \Phi} \lambda(z), \quad z \in G\, ,
\end{equation*}
 induces  a regular
 conformal metric with constant curvature $-4$ on $G$.
 \end{theorem}

We remark that if $\Phi$ is the Perron family of all SK--metrics on $G$, then
$\lambda_{\Phi }(z)\, |dz|$ is the unique maximal conformal metric with
curvature $\le -4$, i.\:\!e.~the hyperbolic metric $\lambda_G(z)\, |dz|$ for
$G$. In particular, it follows that
$\lambda_G(z)\, |dz| \le \lambda_D(z)\,|dz|$ for $z \in D$, if $D\subseteq G$.

\smallskip

\begin{theorem}\label{thm:perron2}
Let $E=(z_j)$ be a sequence of pairwise
distinct points in $G$ and let $(m_j)$ a be sequence of positive integers. Let
\begin{equation*}
 \Phi:= \left\{    \lambda :  \lambda(z)\, |dz| \text{ is an SK--metric on $G$
and } \limsup \limits_{z \to z_j} \frac{\lambda(z)}{|z-z_j|^{m_j}} < +
\infty \text{ for all j } \right\} \, .
\end{equation*}
If $\Phi \not \equiv \emptyset$, then 
\begin{equation*}
\lambda_{\Phi}(z):= \sup_{\lambda \in \Phi} \lambda(z)\, , \quad z \in G\, ,
\end{equation*}
induces a regular conformal pseudometric on $G$ with
$\kappa_{\lambda_{\Phi}}(z)=-4$ for all $z \in G \backslash
E$. Furthermore, $\lambda_{\Phi}(z)\, |dz|$ has a zero of order $m_j$ at $z_j$
for all $j$.
\end{theorem}

Theorem \ref{thm:perron2} guarantees the existence of a
unique maximal conformal pseudometric $\lambda_{\Phi}(z)\, |dz|$ on $G$ with
preassigned zeros. In other words every conformal pseudometric
$\lambda(z)\, |dz|$ on $G$ with curvature bounded above by $-4$ which
vanishes to at least the prescribed order $m_j$ at each $z_j$ of $E$ is
dominated by $\lambda_{\Phi}$, i.\!\;e.~$\lambda\le \lambda_{\Phi}$ in
$G$. This just means $\lambda_{\Phi}(z)\, |dz|$ takes the r\^{o}le of the hyperbolic
metric and Theorem \ref{thm:perron2} is a refined version of the fundamental
theorem respecting zeros.

\medskip 

Since Theorem \ref{thm:perron2} plays an important r\^{o}le later, we include a
proof for convenience of the reader.

\medskip

{\bf Proof of Theorem \ref{thm:perron2}.}
We first note that $\Phi$ is a Perron family of SK--metrics on $G\backslash
E$. Thus $\lambda_{\Phi}(z)$ is well--defined on $G$ and induces a regular conformal 
pseudometric on $G$ with
$\kappa_{\lambda_{\Phi}}(z)=-4$ for $z \in G \backslash E$, see Theorem \ref{thm:perron0}. 

\smallskip

Pick $z_j \in E$ and
choose an open disk $K:=K_{r_j}(z_j)=\{z : |z-z_j|< r_j\}$  such that  $K$ is compactly contained
in $G$ and $K \cap (E\backslash \{z_j\}) = \emptyset$.

\smallskip
To prove that $\lambda_{\Phi}$ has a zero of order $m_j$ at $z_j$ we 
first show that $\lambda_{\Phi} \in \Phi$. 
Note that 
there exists by the fundamental theorem some constant $c$ such that 
$\lambda(z) \le c,\,   z \in K, $
for all $\lambda \in \Phi$. We now define on $K$ the function
\begin{equation*}
\sigma(z):=c\, \left( \frac{|z-z_j|}{r_j} \right)^{m_j}\, .
\end{equation*}
For a fixed $\lambda \in \Phi$ we consider the nonnegative function
$s(z):=\log^+\left(\lambda(z)/\sigma(z)\right)$ on
 $K \backslash \{z_j\}$. Observe that $s$  is subharmonic on
$K\backslash \{z_j\}$ and  since $s$ is bounded at $z_j$ it has a subharmonic
extension to $K$. By construction
$ \limsup _{z \to \zeta} s(z)=0$ for all $\zeta \in \partial K$.
 Hence, $\lambda(z) \le \sigma(z) $ for $z \in K$. As this
holds for every $\lambda \in \Phi$ we obtain
$\lambda_{\Phi}(z) \le \sigma(z)$  
for $z \in K$. Thus we conclude that $\lambda_{\Phi} \in \Phi$. 

\smallskip

By Theorem \ref{thm:existence2} there exists a regular conformal metric
$\mu(z)\, |dz|$ on $K$
 with curvature 
\begin{equation*}
\kappa_{\mu}(z)=-4\, \left( \frac{|z-z_j|}{r_j} \right)^{2m_j}
\end{equation*}
such that $\mu$ is continuous on the closure $\overline{K}$ and 
$\mu\equiv \lambda_{\Phi}$ on $\partial K$.
Thus 
\begin{equation*}
\nu(z):= \left( \frac{|z-z_j|}{r_j} \right)^{m_j}\, \mu(z)
\end{equation*} 
induces a regular conformal pseudometric on $K$ with constant curvature
$-4$, $\nu$ is continuous on $\overline{K}$ and $\nu\equiv \lambda_{\Phi}$ on $ \partial K$.

\smallskip

We are now going to show that $\lambda_{\Phi} (z)= \nu(z)$ for $z \in K$. To
do this, we define $\tilde{s}(z):= \log^+( \lambda_{\Phi}(z)/\nu(z) )$  for $z \in K$. 
Similarly as above, we can conclude that $\tilde{ s}$ is a nonnegative
subharmonic function on $K$. The boundary condition on $\nu$ implies $\lim_{z \to \zeta} \tilde{s}(z)=0$ for
all $\zeta \in \partial K$. So $\lambda_{\Phi}(z) \le \nu(z)$ for $z \in K$.
On the other hand, since $\lambda_{\Phi} \in \Phi$, the gluing lemma (Lemma \ref{lem:gluing}) guarantees that
\begin{equation*}
\tau(z):= 
\begin{cases}
\max\{\lambda_{\Phi}(z), \nu(z) \}\, , & z \in K\, ,\\[2mm]
\lambda_{\Phi}(z)\, ,   &  z \in G\backslash K, 
\end{cases}
\end{equation*}
belongs to $\Phi$. Thus $\nu(z)\le \lambda_{\Phi}(z)$ for $z \in K$. 
Consequently, $\lambda_{\Phi}(z) \, |dz|$ has a zero of order $m_j$ at $z_j$,
which completes the proof.
\hfill{$\blacksquare$}

\medskip

We conclude this section with a result, similar to case of equality in the
fundamental theorem, that is a ``strong version of Ahlfors' lemma'', see \cite[\S
7]{Hei62} and \cite{Chen2001,KR2008,Min1987,Roy1986}.

\begin{lemma}\label{lem:gen_max_3}
Let ${\cal C}$ be a sequence of points in $G$ and let $\lambda(z)\, |dz|$ and $\mu(z)\,
|dz|$ be conformal pseudometrics  
on $G$ with constant curvature $-4$. Suppose ${\cal C}$ is the zero set of
$\mu(z)\, |dz|$ and $\lambda(z)\le
  \mu(z)$ for all $z \in G$.
If
\begin{equation}\label{eq:equal}
\lim_{z\to z_0} \frac{\lambda(z)}{\mu(z)}=1 
\end{equation} 
for some $z_0 \in G$, then $\lambda\equiv \mu$. 
\end{lemma}

\smallskip

{\bf Proof.}
We observe that we can proceed as in the case of ${\cal C}=\emptyset$ 
if (\ref{eq:equal}) is fulfilled for $z_0 \in G
\backslash {\cal C}$ as well as  if (\ref{eq:equal}) is valid for some $z_0\in
{\cal C}$ provided that 
the function $z \mapsto \lambda(z)/\mu(z)$ has a  twice continuously
differentiable extension to a neighborhood of $z_0$.

\smallskip

To see this let $\nu(z)\, |dz|$ be conformal pseudometric  on
some open disk
$K:=K_r(z_0)$  which has constant curvature
$-4$ on $K \backslash \{z_0\}$ and a zero of order $m_0$ at $z_0$. 
 We can further suppose
that $\nu$ is continuous on $\overline{K}$. Now let  $\tilde{\nu}$ be the
continuous extension of $\nu(z)\,|z-z_0|^{-m_0}$ on $\overline {K}$. 
We will show that $\tilde{\nu}$ is twice continuously differentiable on
$K$. For this we note that $\tilde{\nu}$ induces a regular conformal metric on
$K\backslash \{z_0\}$ with curvature
  $\kappa_{\tilde{\nu}}(z)=-4 \, |z-z_0|^{2m_0}$. By Theorem \ref{thm:existence2}
there exists a regular conformal metric $\tau(z)\, |dz|$ on $K$ with
curvature $\kappa_{\tau}(z) = -4\, |z-z_0|^{2m_0}$ which is continuous on
$\overline{K}$ and satisfies
  $\tau(\zeta)=\tilde{\nu}(\zeta)$ for all $\zeta \in \partial K$.
Then the  nonnegative function $s(z):=\log^+(\tilde{\nu}(z)/\tau(z))$  
is  subharmonic first on $K \backslash \{z_0\}$
  and because $s$ is bounded near $z_0$ it extends to a subharmonic function on $K$. Since, by construction,
  $\limsup_{z \to \zeta} s(z)=0$ for all $\zeta \in \partial K$ we
  deduce that $\tilde{\nu}(z) \le \tau(z)$ for all $z \in
  K$. Switching the r\^oles of $\tilde{\nu}$ and $\tau$, we get
  $\tilde{\nu}\equiv \tau$.

\smallskip

Now suppose that (\ref{eq:equal}) holds for some $z_0 \in {\cal
  C}$. If $m_0$ denotes the multiplicity of $z_0$ in ${\cal C}$,  then
(\ref{eq:equal}) implies that $\lambda(z)\, |dz|$ and $\mu(z)\, |dz|$ have a
zero of order $m_0$ at $z_0$.
Thus $\lambda(z)\, |dz|$ and $\mu(z)\, |dz|$ enjoy the same
properties as $\nu(z)\, |dz|$ and the desired result
follows.\hfill{$\blacksquare$}

\section{On the Berger--Nirenberg problem for planar
  domains} 
\label{sec:Berger_Nirenberg_problem}

\subsection{Results}
\label{sec:results}
Suppose $D$ is a regular\footnote{i.\!\;e.~there exists 
  Green's function for $D$ which vanishes
  continuously on $\partial D$.} and bounded domain and $k$ a nonnegative, bounded and (locally) H\"older
continuous function. Under these assumptions it is well--known that there is
always a solution\footnote{A function $u: D \to \R$ is called solution to
(\ref{eq:curvature1}) on $D$, if 
$u \in C^2(D)$ 
and $u$ satisfies (\ref{eq:curvature1}) in $D$.} to the Gauss curvature
equation
\begin{equation}\label{eq:curvature1}
\Delta u= k(z)\, e^{2u}
\end{equation}
on $D$.
On the other hand, if $D$ or $k$ is unbounded, there might be no
solution to (\ref{eq:curvature1}).
For example,
take $D=\C$ and $k(z) = 4 \,|f(z)|^2$ for some entire function
$f\not\equiv 0$. 
By Remark \ref{rem:liouville}, any
solution $u$ would be of the form 
\begin{equation*}
u(z)=\log \left(\frac{1}{|f(z)|}\, \frac{|g'(z)|}{1-|g(z)|^2} \right)\,, \quad z \in \C\,\! ,
\end{equation*}
for some analytic function $g: \C \to \D$. The fact  that a bounded entire
function is constant would then imply that $u \equiv -\infty$,
violating the fact, that $u$ is a solution to (\ref{eq:curvature1}). 

\medskip

In our first result we give for 
regular and bounded domains $D$
 necessary as well as  sufficient conditions on the function
$k$ for the existence of a solution to (\ref{eq:curvature1}) on $D$. 
In the following $g_D$ denotes Green's
function for $D$.

\begin{theorem}\label{thm:sol1}
Let $D$ be a bounded and regular domain and let $k$ be a nonnegative locally H\"older
continuous function on $D$. 
\begin{itemize}
\item[(1)]
If for some (and therefore for every) $z_0 \in D$
\begin{equation*}
\iint \limits_{D} g_D(z_0, \xi)\, k(\xi) \, d\sigma_{\xi} < + \infty\, ,
\end{equation*}
then (\ref{eq:curvature1}) has a solution $u : D \to \R$, which is bounded from above.

\item[(2)]
If (\ref{eq:curvature1}) has a solution $u: D \to \R$ which is bounded from
below and has a harmonic majorant, then  
\begin{equation*}
\iint \limits_{D} g_D(z, \xi) \, k(\xi) \, d\sigma_{\xi} < + \infty
\end{equation*}
for all $z \in D$.
\item[(3)]
There exists a bounded solution 
$u:D \to \R$  to (\ref{eq:curvature1})  if and only if 
\begin{equation*}
\sup \limits_{z \in D} \iint \limits_{D} g_D(z, \xi) \,   k(\xi)\,   d\sigma_{\xi} < + \infty\, .
\end{equation*}
\end{itemize}
\end{theorem}

Now let $D= \D$ in Theorem \ref{thm:sol1}. Then using the elementary estimate
\begin{equation*} 
 \frac{1-|\xi|^2}{2} \le  \log \frac{1}{|\xi|} \le
\frac{1-|\xi|^2}{2\,|\xi|}\, , \quad 0< |\xi| <1,
\end{equation*}
for $g_{\D}(0,\xi)=-\log|\xi|$ 
leads to the following equivalent formulation of Theorem \ref{thm:sol1}.

\begin{corollary}\label{cor:solutions1}
Let $k$ be a nonnegative locally H\"older continuous function on $\D$. 
\begin{itemize}
\item[(1)]
If 
\begin{equation*}
\iint \limits_{\D} (1-|\xi|^2 )\, k(\xi) \, d\sigma_{\xi} < + \infty\, ,
\end{equation*}
then (\ref{eq:curvature1}) has a solution $u : \D \to \R$, which is bounded from above.

\item[(2)]
If (\ref{eq:curvature1}) has a solution $u: \D \to \R$ which is bounded from
below and has a harmonic majorant, then 
\begin{equation*}
\iint \limits_{\D} (1-|\xi|^2 )\, k(\xi) \, d\sigma_{\xi} < + \infty\, .
\end{equation*}
\item[(3)]
There exists a bounded solution $u: \D \to \R$  to (\ref{eq:curvature1}) if and only if 
\begin{equation*}
\sup \limits_{z \in \D} \iint \limits_{\D} \log\left| \frac{1- \overline{\xi} z}{z -
      \xi} \right|\,  k(\xi)\,   d\sigma_{\xi} < + \infty\, .
\end{equation*}
\end{itemize}
\end{corollary}

\bigskip

It might be worth making some remarks on Theorem \ref{thm:sol1} and  Corollary
\ref{cor:solutions1}.
Both, Theorem \ref{thm:sol1} and  Corollary \ref{cor:solutions1}, are not best
possible, because (\ref{eq:curvature1}) may indeed have solutions, even if 
\begin{equation*}
\iint \limits_{D} g_{D}(z, \xi)\, k(\xi) \, d\sigma_{\xi} =+ \infty
\end{equation*}  
for some (and therefore for all) $z \in D$.
Here is an explicit example.

\begin{example} \label{ex:0}
For $\alpha \ge 3/2$ define
\begin{equation*}
\varphi(z)= \frac{1}{(z-1)^{\alpha}}
\end{equation*}
for $z \in \D$ and set $k(z)=4\, |\varphi(z)|^2$ for $z \in \D$.
Then an easy computation yields
\begin{equation*}
\iint \limits_{\D} (1-|z|^2 )\, k(z) \, d\sigma_{z}= +\infty\, .
\end{equation*}
On the other hand, a straightforward check shows that for every analytic and
locally univalent self--map $f$ of $\D$  the function
\begin{equation*}
u_f(z):=  \log \left( \frac{1}{|\varphi(z)|} \, \, \frac{|f'(z)|}{1-|f(z)|^2}
\right) 
\end{equation*}
is a solution to
(\ref{eq:curvature1}) on $\D$.
\end{example}

Observe that in Example \ref{ex:0} the  function $k$ 
is the squared modulus of a holomorphic function.
Thus Theorem \ref{thm:d} applies and shows that (\ref{eq:curvature1}) does have
solutions on $\D$.

\bigskip

A second remark is that Theorem \ref{thm:sol1} (c)  characterizes those functions $k$ 
for which (\ref{eq:curvature1}) has at least one
bounded solution. In particular, when $D=\D$ and $k=4\,|f|^2$ for some
holomorphic function  $f$ in $\D$, then we have the
following connection.

\begin{remark}\label{rem:boundedsol}
Let $\varphi: \D \to \C$  be analytic and $k(z)=4 \, |\varphi'(z)|^2$.
Then there exists a bounded solution to (\ref{eq:curvature1}) if and only if
$\varphi \in BMOA$, where 
\begin{equation*}
BMOA=\left\{ \varphi: \D \to \C \text{ analytic} \, : \,  \sup_{ z \in \D} \iint \limits_{\D} g_{\D}(z,
  \xi)\, |\varphi'(\xi)|^2 \, d\sigma_{\xi} < + \infty \right \}
\end{equation*}
is the space of analytic functions of bounded
mean oscillation on $\D$, see \cite[p.~314/315]{BF2008}.
\end{remark}

We further note that in Theorem \ref{thm:sol1} and Corollary
\ref{cor:solutions1} condition (1) does not imply condition (3). 
The Gauss curvature equation (\ref{eq:curvature1})   may indeed
have solutions but none of the solutions is bounded. For example choose
$\varphi \in H^2\backslash BMOA$,
where $H^2$ denotes the Hardy space consisting of the functions $f$
analytic in $\D$ for which  the integrals  
\begin{equation*}
\int \limits_0^{2\pi} |f(re^{it})|^2\, dt  
\end{equation*}
remain bounded as $r \to 1$.
By Littlewood--Paley's identity, cf.~Remark \ref{rem:little_paley}, it follows that $\varphi'
\in {\cal A}_1^2$. Now set $k(z)=4\, |\varphi'(z)|^2$. Then by Theorem
\ref{thm:d} the Gauss curvature equation (\ref{eq:curvature1}) does have
solutions and according to Remark \ref{rem:boundedsol} every solution to
(\ref{eq:curvature1}) must be unbounded.

\medskip

Finally, suppose $k$ is a nonnegative, locally H\"older continuous and radially symmetric function on $\D$. Then 
 Corollary \ref{cor:solutions1} allows us to characterize those functions $k$ for which the
 Gauss curvature equation (\ref{eq:curvature1}) has a solution on $\D$ with a harmonic majorant.

\begin{corollary}\label{cor:radialsol}
Let $k$ be a nonnegative locally H\"older continuous function on $\D$ such that
$k(\xi)=k(|\xi|)$ for all $\xi \in \D$. Then (\ref{eq:curvature1}) has a solution $u:\D \to \R$  with a harmonic majorant  if and only if
 \begin{equation}\label{eq:k_radsym}
\iint \limits_{\D} (1-|\xi|^2 )\, k(\xi) \, d\sigma_{\xi} < + \infty\, .
\end{equation}
\end{corollary}

\medskip

To illustrate the use of Corollary \ref{cor:solutions1} and Corollary
\ref{cor:radialsol} respectively, here is an example.

\begin{example}\label{ex:k1}
Let $\gamma \in \R$, $\gamma \ge 1$, and define for $z \in \D$ 
\begin{alignat*}{1}
k_1^{\gamma}(z)= & \frac{1}{(1-|z|^2)^2}\,  \frac{1}{\left[\log \left( \frac{e}{1-|z|^2}
  \right)\right]^{\gamma}}\\[2mm]
k_2^{\gamma}(z) =  & \frac{1}{(1-|z|^2)^2}\,  \frac{1}{\log \left( \frac{e}{1-|z|^2}
  \right)}\, \frac{1}{\left[\log \left(e \log \left( \frac{e}{1-|z|^2}
  \right)\right)\right]^{\gamma}}\\[2mm]
k_3^{\gamma}(z)= &\frac{1}{(1-|z|^2)^2}\,  \frac{1}{\log \left( \frac{e}{1-|z|^2}
  \right)}\, \frac{1}{\log \left(e \log \left( \frac{e}{1-|z|^2}
  \right)\right)}\, \frac{1}{\left[\log\left( e \log \left(e \log \left( \frac{e}{1-|z|^2}
  \right)\right)\right)\right]^{\gamma}}\\[2mm]
\text{etc.\!\;.} \quad \, \,&
\end{alignat*}

\begin{itemize}
\item[(a)]
If  a nonnegative locally H\"older continuous function $k$ on $\D$  satisfies
$k(z) \le k_j^{\gamma}(z)$ for all $z \in \D$ and some $\gamma > 1$ and $j \in \N$, then the
Gauss curvature equation (\ref{eq:curvature1}) has a solution on $\D$. 

\item[(b)]
If $k$ is a continuous function on $\D$ such that  $k(z) \ge k_j^1(z)$ for all $z
\in \D$ and some
$j \in \N$, then there is no solution to (\ref{eq:curvature1}) on $\D$.
\end{itemize}
\end{example}

\smallskip

A few comments are in order. First, part (a) of Example \ref{ex:k1} is a
consequence of
 Corollary \ref{cor:solutions1} (1). In fact, a straightforward computation
gives
\begin{equation*}
\iint \limits_{\D} (1-|z|^2)\, k_j^{\gamma}(z)\, d\sigma_z < +\infty
\end{equation*}
for every $\gamma >1$ and $j \in \N$. 

\smallskip

For the special case that $k$ is an essentially positive\footnote{A locally
  H\"older continuous function is called essentially positive, if there is a
  strictly increasing sequence $(G_n)$ of relatively compact subdomains $G_n$ of $\D$ such that
  $\D=\cup_{n} G_n$ and $k(\zeta)>0$ for $\zeta \in \partial G_n$ for all
  $n\in \N$.}\label{page:ess} function
Example \ref{ex:k1} (a) has been discussed by Kalka and Yang 
\cite[Theorem 3.1]{KY1993}.
This additional hypothesis on $k$ guarantees a solution to
(\ref{eq:curvature1}) on $\D$ which even tends to $+ \infty$ at the boundary of $\D$.
For the proof Kalka and Yang  use a generalized Perron method. In particular, the existence of a
subsolution\footnote{A  function $u: G \to \R$ is said to be a
subsolution to (\ref{eq:curvature1}) on $G$,  if  $u \in C^2(G)$ and
$\Delta u \ge
k(z) \, e^{2\!\; u}$ on $G$.} to the Gauss curvature equation (\ref{eq:curvature1}) on $\D$ guarantees a
solution to (\ref{eq:curvature1}) on $\D$. 
The authors give for each $\gamma >1$ and $j \in \N$ an explicit subsolution
to the PDE $\Delta u =k_j^{\gamma}(z)\, e^{2u}$ on $\D$.

\smallskip
Second, statement (b) of Example \ref{ex:k1} is due to  Kalka and Yang \cite[Theorem
3.1]{KY1993}. The key step in their proof  consists in showing that for each
$j \in \N$
there is no 
solution to the Gauss curvature equation $\Delta u =k^1_{j}(z)\, e^{2u}$ on
$\D$. The assertion of Example \ref{ex:k1} (b) follows then directly by  
employing their generalized Perron method. For the key step, 
Kalka and Yang give a  quite intricate argument, which 
relies heavily on Yau's celebrated maximum principle for
complete metrics \cite{Yau1975, Yau1978}. A much simpler, almost elementary
 proof that  there is no solution to (\ref{eq:curvature1}) on $\D$ for
 $k(z)=k^1_{j}(z)$, $j \in \N$, can be found in \cite{Kra2011b}.
This approach has also other ramifications  which are discussed in \cite{Kra2011b}.

\smallskip

Third, we note that  conditions (a) and (b) in Example \ref{ex:k1} do not complement each other.
For example, choose $k(z)=|z+1|^{-2} + |z-1|^{-2}$ for $z \in \D$. Since
$k$ fulfills the hypothesis of Corollary \ref{cor:solutions1} we conclude
that  (\ref{eq:curvature1}) has a solution on $\D$, but condition (a) of
Example \ref{ex:k1} is not applicable. Hence Theorem \ref{thm:sol1},
Corollary \ref{cor:solutions1} and Corollary \ref{cor:radialsol} generalize
the results of Kalka and Yang.

\smallskip

Finally, an easy computation shows that 
\begin{equation*}
\iint \limits_{\D} (1-|\xi|^2 )\, k^1_j(\xi) \, d\sigma_{\xi} = + \infty\, 
\end{equation*}
for $j \in \N$. On the other hand, by Example \ref{ex:k1} (b) 
no solution  to the Gauss
curvature equation (\ref{eq:curvature1})  on $\D$  is possible
for $k(z)=k_j^{1}(z)$, $j \in \N$. Thus one is
inclined to ask whether for a radially symmetric, locally H\"older continuous function
  $k: \D \to [0, +\infty)$ equation (\ref{eq:curvature1}) has no solution on $\D$ 
 if and only if 
 \begin{equation*}
\iint \limits_{\D} (1-|\xi|^2 )\, k(\xi) \, d\sigma_{\xi} = + \infty\, .
\end{equation*}


\subsection{Proof of  Theorem \ref{thm:sol1}}
\label{sec:proofs1}
To prove Theorem \ref{thm:sol1} we need two results. The first is the
solvability of the Dirichlet problem for the Gauss curvature equation
(\ref{eq:curvature1}) and the second is a Harnack type theorem for solutions to
(\ref{eq:curvature1}).

\begin{theorem}\label{thm:existence2}
Let $D$ be a bounded and regular domain, let $k$ be a bounded and nonnegative locally
H\"older continuous function on $D$  and let
$\tau: \partial D \to \R$ be a continuous function.
\begin{itemize}
\item[(a)]
There exists a (unique) function $u \in C(\overline{D}) \cap C^2(D)$ which solves
the boundary value problem
 \begin{equation}\label{eq:ex1}
\begin{array}{rclll}
\Delta u&=&k(z) \, e^{2 u} &\text{in  } &  D\!\;,\\[2mm]
  u &\equiv & \tau    & \text{on  } & \partial D\!\;.
\end{array}
\end{equation}
In particular,
\begin{equation}\label{eq:ex2}
u(z)=h(z)- \frac{1}{2\pi}\iint \limits_{D} g_D(z, \xi)\, k(\xi)\,
e^{2\!\;u(\xi)}\, d\sigma_{\xi}\, , \quad z \in D\, ,
\end{equation}
where $h$ is harmonic in $D$ and continuous on $\overline{D}$ satisfying $h \equiv
\tau$ on $\partial D$.
\item[(b)]
If, conversely, a bounded and integrable function $u$ on $D$ satisfies
(\ref{eq:ex2}), then $u$ belongs to $C(\overline{D}) \cap C^2(D)$ and solves (\ref{eq:ex1}).
\end{itemize}
\end{theorem}
 
For a proof of Theorem \ref{thm:existence2} we refer the reader to
\cite[p.~286]{Cou1968} and
\cite[p.~53--55 $\&$ p.~304]{GT}.

\medskip

\begin{lemma} \label{cor:monosequence}
Let  $k$ 
be a nonnegative locally H\"older continuous function on $G$ and
$(u_n)$ be a monotonically decreasing  sequence of solutions to
(\ref{eq:curvature1}) on $G$. If $\lim_{n \to \infty} u_n(z_0) =-\infty$ for some $z_0
\in G$, then $(u_n)$
converges locally uniformly in $G$ to $-\infty$, otherwise $(u_n)$
converges locally uniformly in $G$ to $u:=\lim_{n \to \infty} u_n$
and $u$ is a solution to (\ref{eq:curvature1}) on $G$.
\end{lemma}

A proof of Lemma \ref{cor:monosequence} for the special case $k \equiv 4$ can be found in
\cite[\S 11]{Hei62}. The proof for the more general situation needs only
slight modifications and will therefore be omitted. See also \cite[Proposition 4.1]{MT2002}.

\medskip

{\bf Proof of Theorem \ref{thm:sol1}.}\\
(1) Let $(D_n)$ be a sequence of relatively compact regular subdomains of $D$ such that
$z_0 \in D_1\subset D_2 \subset D_3 \subset \ldots$ and $D =\cup_n \,
D_n$. Pick a constant $c \in \R$. 
 For each $n$, let $u_n: D_n\to
  \R$ denote the solution to the boundary value problem
\[
\begin{array}{rclll}
\Delta u&=&k(z) \, e^{2 u} &\text{in  } &  D_n\, ,\\[2mm]
  u &\equiv & c    & \text{on  } & \partial D_n\, .
\end{array}
\]
Thus we can write
\begin{equation*}
u_n(z)=c-\frac{1}{2 \pi} \iint \limits_{D_n}g_{D_n}(z, \xi)\, k(\xi) \,
e^{2u_n(\xi)} \, d\sigma_{\xi}\, , \quad \, z \in D_n\, .
\end{equation*}
Since $u_n$ is subharmonic on $D_n$ and $g_{D_n}(z,
\xi)\le g_D(z, \xi)$ for all $z, \xi \in D_n$, we obtain
\begin{equation}\label{eq:1}
\begin{split}
 u_n(z_0)&=c-\frac{1}{2 \pi} \iint \limits_{D_n} g_{D_n}(z_0, \xi) \, k(\xi) \,
e^{2u_n(\xi)} \, d\sigma_{\xi}\\[2mm]
&\ge c- e^{2c } \, \iint \limits_{D} g_D(z_0,
\xi) \, k(\xi) \,  d\sigma_{\xi} \ge \tilde{c}
\end{split}
\end{equation}
for some finite constant $\tilde{c}$. Letting $n \to \infty$ yields 
\begin{equation*}
\liminf\limits_{n \to \infty} u_n(z_0) > - \infty\, .
\end{equation*}
Note that the boundary condition on $u_n$ implies that $(u_n)$ is a monotonically decreasing sequence of solutions to $\Delta
u=k(z)\, e^{2\!\;u}$.
 Thus Lemma
\ref{cor:monosequence} applies and 
\begin{equation*}
u(z):= \lim_{n \to \infty} u_n(z)\, , \quad z \in D\,,
\end{equation*}   
is a solution to (\ref{eq:curvature1}) on $D$, which is 
bounded above by construction. 

\smallskip

(2) Let $u : D \to \R$ be a solution to (\ref{eq:curvature1}) with some
harmonic majorant. Then by  the Poisson--Jensen formula
\begin{equation*} 
u(z)=h(z)-\frac{1}{2 \pi} \iint \limits_{D}g_D(z, \xi)\, k(\xi)\,
e^{2u(\xi)} \, d\sigma_{\xi}
\end{equation*}
for $z \in D$, where $h$ is the least harmonic majorant of $u$ on $D$. As
$u(z)>c> -\infty $ for $z \in D$, we get for  fixed $z \in D$
\begin{equation*} 
\begin{split}
\frac{1}{2 \pi} \iint \limits_{D}g_D(z, \xi)\, k(\xi)
\, d\sigma_{\xi} &\le e^{-2 c} \, \frac{1}{2 \pi} \iint \limits_{D}g_D(z, \xi)\, k(\xi)\,
e^{2 u(\xi)} \, d\sigma_{\xi}\\[2mm] &= e^{-2 c}\, \big(h(z)-u(z)\big
) < + \infty\,
,
\end{split}
\end{equation*}
as desired.

\smallskip

(3) Let $u: D \to \R$ be a bounded solution to (\ref{eq:curvature1}),
i.\!\;e.~$|u(z)| \le c$ for $z \in D$, where $c$ is 
some positive constant. Hence
\begin{equation*}
\begin{split}
\frac{1}{2 \pi} \iint \limits_{D} g_D(z, \xi) \, k(\xi)\, d\sigma_{\xi} &\le   
e^{2 c}\, \frac{1}{2 \pi} \iint \limits_{D} g_D(z, \xi) \, k(\xi)\, e^{2 u(\xi)}\,
  d\sigma_{\xi}\\[2mm]& =  e^{2c} \, \big( h(z) -u(z)\big) \le 2 c\,
  e^{2c}
\end{split}
\end{equation*}
for $z \in D$,
where $h$ is the least harmonic majorant of $u$ on $D$.
\smallskip

Conversely, suppose that 
\begin{equation}\label{eq:char1}
\sup \limits_{z \in D} \iint \limits_{D} g_D(z, \xi) \,   k(\xi)\,   d\sigma_{\xi} < + \infty\, .
\end{equation}
Let $(u_n)$ be the sequence  constructed in (1). Then 
\begin{equation*}
u(z):=\lim \limits_{n \to \infty} u_n(z)\, , \quad z \in D,
\end{equation*}
is a solution to (\ref{eq:curvature1}), which is bounded from above. Inequality (\ref{eq:1})
combined with (\ref{eq:char1}) shows that there is a constant $c_1$ such that
$u_n(z) \ge c_1$ for all  $ z \in D$ and all
$n$, so $u$ is also bounded from below.
\hfill{$\blacksquare$}

\subsection{Proof of Theorem \ref{thm:d} and Corollary \ref{cor:radialsol}}
\label{sec:proofs2}

The proof of Theorem \ref{thm:d} relies on Liouville's theorem (Remark
\ref{rem:liouville}) and the Littlewood--Paley identity.
Before beginning the proof, it will be useful to recall the Littlewood--Paley identity
\cite[p.~178]{Sha1993}.

\begin{remark}[Littlewood--Paley's identity]\label{rem:little_paley}
Let $\varphi: \D \to \C$ be a holomorphic function. Then we have
\begin{equation*}
\frac{1}{2\pi} \int \limits_0^{2\pi} |\varphi(e^{it})|^2\, dt= |\varphi(0)|^2 +
\frac{2}{\pi} \iint \limits_{\D} \log \frac{1}{|z|} \, |\varphi'(z)|^2\,
d\sigma_{z}\, .
\end{equation*}
 This in particular shows that
\begin{equation*}
{\cal A}_1^2=\{\varphi': \varphi \in H^2   \}\, .
\end{equation*}

\end{remark}

\bigskip

{\bf Proof of Theorem \ref{thm:d}.}
We will show that for an analytic function $\varphi$ on $\D$ the Gauss
curvature equation
\begin{equation}\label{eq:curvature_holomorph}
\Delta u = 4\, |\varphi(z)|^2 e^{2u}
\end{equation}  
has a solution $u:\D \to \R$ if and only if 
$\varphi(z)=\varphi_1(z)\, \varphi_2(z)$ for some function $\varphi_1 \in
{\cal A}_1^2$  and a nonvanishing analytic function $\varphi_2: \D
 \to \C$.

\medskip

We first
note that it suffices to consider the case $\varphi \not \equiv 0$. Now
suppose $\varphi \not \equiv 0$ and 
 $u: \D \to \R$ is a solution to (\ref{eq:curvature_holomorph}). Then, by
Liouville's theorem,
\begin{equation*}
u(z)=\log \left( \frac{1}{|\varphi(z)|}\, \frac{|f'(z)|}{1-|f(z)|^2} \right)\,,  \quad z
\in \D\, ,
\end{equation*}
for some analytic self--map $f$ of $\D$.
Since $\varphi$ and $f'$ have the same zeros, we can write $\varphi(z)=f'(z)\,
\varphi_2(z)$, where  $\varphi_2$ is analytic and zerofree in $\D$. 
From Remark \ref{rem:little_paley} it follows that $f' \in {\cal A}_1^2$.

\smallskip

Conversely, let $\varphi(z)=\varphi_1(z)\, \varphi_2(z)$, where $\varphi_1 \in
{\cal A}_1^2$, $\varphi_1 \not \equiv 0$, and $\varphi_2: \D \to \C\backslash \{ 0 \}$ is an analytic function. Then 
Corollary \ref{cor:solutions1} (1) ensures a solution $u$ to
$\Delta u= 4\, |\varphi_1(z)|^2\, e^{2\!\;u}$ on $\D$. By Liouville's theorem
there exists an analytic function $f:\D \to \D$ such that
\begin{equation*}
u(z)=\log \left(\frac{1}{|\varphi_1(z)|}\, \frac{|f'(z)|}{1-|f(z)|^2}\right)\,,  \quad z
\in \D\, .
\end{equation*}
Hence the function $f'/\varphi_1$ is analytic and zerofree in $\D$. So,
\begin{equation*}
\tilde{u}(z):= \log \left(\frac{1}{|\varphi(z)|}\,
  \frac{|f'(z)|}{1-|f(z)|^2} \right)\,,
\quad z \in \D\, ,
\end{equation*}
is well--defined and  a solution to $\Delta u= 4\, |\varphi(z)|^2\,
e^{2\!\;u}$ on $\D$.\hfill{$\blacksquare$}

\medskip
 
If $k$ is a nonnegative, locally H\"older continuous and  radially symmetric
function on $\D$, then only the existence of a solution to equation
(\ref{eq:curvature1}) on $\D$ with a harmonic majorant yields (\ref{eq:k_radsym}). 
Thus, in this special case the hypothesis in part  (2) of Corollary \ref{cor:solutions1} can slightly be relaxed. 

\medskip 

{\bf Proof of Corollary \ref{cor:radialsol}.} By Corollary \ref{cor:solutions1} (1) condition (\ref{eq:k_radsym}) ensures a
solution to (\ref{eq:curvature1}) with a harmonic majorant. 

\smallskip

For the converse, let $u$ be a solution to (\ref{eq:curvature1}) on $\D$ with a
harmonic majorant. Then it follows by Green's theorem and Jensen's inequality
that the function 
\begin{equation*}
v(z):=\frac{1}{2\pi} \int \limits_{0}^{2\pi} u(|z|e^{it})\, dt\, ,\quad  z \in
\D\, ,
\end{equation*}
is a subsolution to (\ref{eq:curvature1}), i.\!\;e.~$v \in C^2(\D)$ and
satisfies $\Delta v \ge k(z) \, e^{2v}$ on $\D$. Note, that $v$ has a harmonic majorant,
see \cite[Chapter I, Theorem 6.7]{Gar2007}. As
$v$ is subharmonic in $\D$,  we  apply  the Poisson--Jensen formula and  get
with the help of Jensen's inequality
\begin{equation*}
\begin{split}
h(0)-v(0)&= \frac{1}{2\pi} \iint \limits_{\D}\log\frac{1}{|\xi|}\,  \Delta
v(\xi) \, d \sigma _{\xi}
\ge \frac{1}{2\pi} \iint \limits_{\D} \log\frac{1}{|\xi|}\, k(\xi)\,  e^{2\!\;v(\xi)} \, d
\sigma _{\xi}\\[2mm]
& \ge e^{2\!\;v(0)} \, \frac{1}{2\pi} \iint \limits_{\D} \log\frac{1}{|\xi|}\, k(\xi) \, d
\sigma _{\xi}\, ,
\end{split}
\end{equation*}
where $h$ is the least harmonic majorant of $v$ on $\D$. Hence
(\ref{eq:k_radsym}) holds. \hfill{$\blacksquare$}

\section{Maximal conformal  pseudometrics and maximal functions}
\label{sec:maximal}

Let ${\cal C}=(z_1, \ldots, z_1, z_2, \ldots, z_2, \ldots)$, $z_j
\not=z_k$ for $j\not=k$, be a sequence in $\D$  and denote by $m_j$  the
multiplicity of $z_j$ in ${\cal C}$. Suppose that the family $\Phi_{\cal C}$ of
all SK--metrics $\lambda(z)\, |dz|$ which vanishes at least on ${\cal C}$,
i.\!\;e.~
\begin{equation*}
\limsup_{z\to z_j} \frac{\lambda(z)}{|z-z_j|^{m_j}}<+ \infty
\end{equation*}
for all $j$,
is not empty.
Then, by Theorem \ref{thm:perron2},
\begin{equation*}
\lambda_{max}(z)\, |dz|:= \sup_{\lambda \in \Phi_{\cal C} } \lambda(z) \,|dz|\,
, \quad z \in \D\, ,
\end{equation*}
defines the unique maximal conformal pseudometric on $\D$ with constant curvature
 $-4$ and zero set ${\cal C}$.
Even though   the existence of
 $\lambda_{max}(z)\, |dz|$ is guaranteed if ${\cal C}$ is the  zero set of an ${\cal
   A}_1^2$ function, see  Theorem \ref{thm:main},  we lack in explicit examples. 
One way round this problem is Liouville's theorem. In particular, we are
interested in 
the developing map $F$ of 
a maximal conformal metric $\lambda_{max}(z)\, |dz|$, i.\!\;e. 
\begin{equation*}
\lambda_{max}(z)= \frac{|F'(z)|}{1-|F(z)|^2}\, ,\quad z \in \D\, ,
\end{equation*}
for some analytic function $F: \D \to \D$.
The analytic functions which represent  
maximal conformal metrics are of some interest
in their own right as they are natural generalizations of the unit disk
automorphisms, i.\!\;e.~the developing maps of the Poincar\'e metric $\lambda_{\D}(z)\,
|dz|$. Thus the following definition might be appropriate.

\begin{definition}\label{def:maximal-function}
Let ${\cal C}$ be a sequence in $\D$ and assume that
 $\lambda_{max}(z)\, |dz|$ is the maximal conformal pseudometric for $\D$ 
with constant curvature $-4$ and zero set ${\cal C}$. 
Then every developing map $F$ of $\lambda_{max}(z)\, |dz|$
is called maximal function with critical set ${\cal C}$.
\end{definition}

Note that by Liouville's theorem a maximal function is uniquely determined by
its critical set 
 up to 
postcomposition with a unit disk automorphism.

\medskip

As an immediate consequence of Theorem \ref{thm:0}, i.\!\;e.~every maximal function  is an indestructible Blaschke product, we obtain the equivalence of 
statements (a) and (b) in Theorem \ref{thm:main}, which we now
restate for convenience of reference.

\medskip
\begin{corollary} \label{cor:1}
 Let ${\cal C}$ be a sequence of points in $\D$. Then the following are
equivalent.
\begin{itemize}
\item[(a)] There exists an analytic self--map of $\D$ with critical set
  ${\cal C}$.
\item[(b)] There exists an indestructible Blaschke product with critical set
  ${\cal C}$.
\end{itemize}
\end{corollary}

{\bf Proof.} \\
(a) $\Rightarrow$ (b):
 Let $f :\D \to \D$ be analytic with critical set
 ${\cal C}$. Then
\begin{equation*}
\lambda(z)\, |dz|= \frac{|f'(z)|}{1-|f(z)|^2}\, |dz|
\end{equation*}
is a conformal pseudometric on $\D$ with constant curvature $-4$ and zero set ${\cal
  C}$. So $\lambda \in \Phi_{\cal{C}}$\footnote{Recall
  $\Phi_{\cal C}$ denotes the family of all SK--metrics which vanishes at
least on ${\cal C}$.} and $\Phi_{\cal C}$ is not empty.
This legitimizes the use of Theorem \ref{thm:perron2} which gives a maximal conformal pseudometric on $\D$ with constant curvature $-4$ and zero set ${\cal
  C}$. Thus there is a maximal 
 function $F:\D \to \D$ with critical set ${\cal C}$, which is an
 indestructible  Blaschke
 product by Theorem \ref{thm:0}. \hfill{$\blacksquare$}

\medskip

\label{page:2}
Theorem \ref{thm:0} and Corollary \ref{cor:1} 
merit some comment.
First, if  ${\cal C} \subset \D$ is a finite
sequence, then Heins observed that the maximal functions for ${\cal C}$ are
precisely the finite Blaschke products with critical set ${\cal C}$, cf.~\cite[\S 29]{Hei62}.
Heins' proof is purely topological and splits into two parts. In the first
he shows that an analytic self--map of $\D$ which has constant finite valence,
i.\!\;e.~a finite
Blaschke product, is a maximal function.
Secondly he establishes the following theorem:  

\smallskip

\begin{satz}\label{thm:finite_blaschke}
Let ${\cal C}$ 
 be a finite sequence in $\D$ that contains $n$ points. 
Then there exists a finite Blaschke product $F$ with critical set
${\cal C}$. $F$ is unique up to
postcomposition with a unit disk automorphism. In this case $F$ has degree
$m=n+1$.
\end{satz}

\smallskip 

The uniqueness statement in Theorem \ref{thm:finite_blaschke} 
follows easily
from Nehari's generalization of Schwarz' lemma \cite[Corollary to Theorem 1]{Neh1946}. 
To settle the existence part Heins showed that the set of critical points of all
finite Blaschke products of degree $n+1$, which is clearly closed, is also open in the
poly disk $\D^n$ by applying Brouwer's fixed point theorem.
Similar  proofs to Theorem \ref{thm:finite_blaschke} 
can also be found in papers by  Wang \& Peng \cite{WP79} and Zakeri
 \cite{Z96}. A completely different approach to Theorem \ref{thm:finite_blaschke}
 via Circle Packing is due to Stephenson, see  \cite[Lemma 13.7 and Theorem 21.1]{Ste2005}. Stephenson builds 
 discrete finite Blaschke products with prescribed branch set and shows that
 under refinement  these
 discrete Blaschke products converge locally uniformly in $\D$ to a classical
 Blaschke product with the desired  critical points.

\medskip

\label{page:3}
A further remark is that if ${\cal C}$ is finite, then  Corollary \ref{cor:1} does not directly imply the
existence statement of Theorem \ref{thm:finite_blaschke}. However, as a consequence of Theorem \ref{thm:0}
we can deduce that if ${\cal C}$ is finite, then  every maximal function for
${\cal C}$ is a
finite Blaschke product with critical set ${\cal C}$, see Theorem
\ref{thm:finiteC} (a).  Now applying  Nehari's uniqueness result we also
arrive at  Heins'
characterization of maximal functions for  finite sequences  ${\cal C}$, but
in  a completely different way.

\medskip
  
Finally, Corollary \ref{cor:1} is discussed in \cite[Theorem 2.1]{KR} for the case that ${\cal C}
\subset \D$ is a Blaschke sequence. There the idea of the proof is the following.
In a first step the existence of a solution $u: \D \to \R$ to the
boundary value problem
\begin{alignat*}{2}
\Delta u &= |B(z)|^2\, e^{2u} & \quad &\text{ in } \D\, ,\\
      \lim_{z \to \zeta} u(z)&=+ \infty & \quad &\text{ for every } \zeta
      \in \partial \D\, ,
\end{alignat*}
is guaranteed,
where $B$ is a Blaschke product whose {\it zero set} consists precisely of the set
${\cal C}$. Then as a consequence of Liouville's theorem there is an analytic self--map $f$ of $\D$
with critical set ${\cal C}$ which represents the solution $u$. The boundary
condition on $u$ then 
implies that $f$ is an inner function and a result by Frostman
(cf.~\cite[Chapter II, Theorem 6.4]{Gar2007}) gives the desired
Blaschke product with critical set ${\cal C}$.

\smallskip

All methods which were employed to prove Theorem \ref{thm:0} and Corollary
\ref{cor:1} for finite and Blaschke sequences ${\cal C}$ seem too restrictive to be
adaptable to the general situation, since they heavily rely on the special
choice of ${\cal C}$.
Our approach to Theorem \ref{thm:0} is exclusively based on conformal
pseudometrics with  curvature bounded above by $-4$, i.\!\;e.~SK--metrics.

\medskip

{\bf Proof of Theorem \ref{thm:0}.}
Let $F: \D \to \D$ be a maximal function with critical set 
\begin{equation*}
{\cal C}=(\underbrace{z_1, \ldots, z_1}_{m_1
  -\text{times}},\underbrace{z_2, \ldots, z_2}_{m_2
  -\text{times}} , \ldots  ) 
\end{equation*}
and let
\begin{equation*}
\lambda_{max}(z) \, |dz|= \frac{|F'(z)|}{1-|F(z)|^2} \, |dz|\, ,  \quad z
\in \D\, ,
\end{equation*}
be the maximal conformal pseudometric of constant
curvature $-4$ with zero set ${\cal C}$. In particular, if $\lambda(z)\,
|dz|$ is a conformal pseudometric on $\D$ with curvature bounded above by $-4$ which satisfies
\begin{equation*}
\limsup_{z \to z_j} \frac{\lambda(z)}{|z-z_j|^{m_j}} < + \infty
\end{equation*}
for all $j$, then $\lambda(z)\le \lambda_{max}(z)$ for $z \in \D$.

\smallskip

To show  that $F$ is a Blaschke product, we use its canonical
factorization, i.\!\;e.~$F=B\, S\, O$, where $B$ is a Blaschke product, $S$
a singular function and $O$ an outer function, see \cite[Chapter II, Corollary 5.7]{Gar2007}. 

\smallskip

We assume that $S\not \equiv \eta$, where $|\eta|=1$. Then
by a result due to Frostman \cite[Chapter II, Theorem 6.2]{Gar2007} there is some $\zeta
\in \partial \D$  such that the angular limit
\begin{equation*}
\angle \lim_{z \to \zeta} S(z)=0\, .
\end{equation*}
In particular, $\angle \lim_{z \to \zeta} F(z)=0$.
Choose $ \alpha \in (-1,0)$ and let
\begin{equation*}
\lambda_{\alpha}(w)\, |dw|:= \frac{(\alpha +1)}{|w|^{|\alpha|}} \, \frac{ 1}{1-|w|^{2\,(\alpha +1)}}\,
|dw|\, .
\end{equation*}
We note that $\lambda_{\alpha}(w)\, |dw|$ is a conformal metric of constant curvature $-4$ on $\D':=\D \backslash \{ 0 \}$.
Thus the pullback of
$\lambda_{\alpha}(z)\, |dz|$ via $F$ restricted to $\D \backslash {\cal C}$ defines 
a conformal metric on $\D\backslash {\cal C}$ with constant curvature $-4$,
 i.\!\;e.
\begin{equation*}
\nu(z)\, |dz|:=\frac{(\alpha +1)}{|F(z)|^{|\alpha|}}\,
\frac{|F'(z)|}{1-|F(z)|^{2( \alpha +1)}}\, |dz|\, .
\end{equation*}
By Lemma \ref{lem:new_sk} we conclude that 
\begin{equation*}
\begin{split}
\sigma(z)\, |dz|&:=\nu(z) \, |B(z)|^{|\alpha|} \,|O(z)|^{|\alpha|}\,|dz|\\[2mm]
& =\frac{(\alpha +1)}{|S(z)|^{|\alpha|}}\, 
\frac{|F'(z)|}{1-|F(z)|^{2( \alpha +1)}} \, |dz|
\end{split}
\end{equation*}
is a conformal metric on $\D \backslash {\cal C}$ with curvature bounded above
by $-4$. Further, by construction $\sigma$ is a
continuous function on $\D$.
Now Lemma \ref{lem:hebsing} guarantees that $\sigma(z)\, |dz|$ is a conformal
pseudometric on $\D$ with curvature bounded above by $-4$. Obviously, ${\cal
  C}$ is the zero set of $\sigma(z)\, |dz|$.
This immediately implies that
\begin{alignat*}{2}
\sigma (z) &\le \lambda_{max}(z) &\text{for } z\in \D \phantom{\, .}
\intertext{ and therefore}
\frac{(\alpha +1)}{|S(z)|^{|\alpha|}}\,
\frac{1}{1-|F(z)|^{2( \alpha +1)}} &\le \frac{1}{1-|F(z)|^{2}}\,
& \quad \text{for } z\in \D\, .
\end{alignat*}
Now
 letting $z \stackrel{\angle}{\to} \zeta$, we get
\begin{equation*}
+\infty\stackrel{z \stackrel{\angle}{\to} \zeta}{\longleftarrow} \frac{(\alpha +1)}{|S(z)|^{|\alpha|}}\,
\frac{1}{1-|F(z)|^{2( \alpha +1)}} \le \frac{1}{1-|F(z)|^{2}}\stackrel{z
  \stackrel{\angle}{\to} \zeta}{\longrightarrow} 1 \, ,
\end{equation*}
the desired contradiction. Hence $F(z)=\eta \, B(z)\, O(z)$ for some constant $\eta$
with $|\eta|=1$.

\smallskip

Since $B$ and $O$ are bounded analytic functions, 
there exists by Fatou's theorem a subset $A \subset \partial \D$ of  full Lebesgue measure
such that
\begin{equation*}
\angle \lim_{ z \to \zeta} B(z) \in \partial \D \quad \text{ and }  \quad
\angle \lim_{ z \to \zeta} O(z) \in \overline{\D}
\end{equation*}
exist for
all $\zeta \in A$. Pick  a point $\zeta \in A$ and assume  $\angle \lim_{ z
  \to \zeta} O(z)= \beta$ where $0\le |\beta| <1$. Thus $\angle \lim_{ z
  \to \zeta} |F(z)|= |\beta| $.
We now consider  the hyperbolic metric $\lambda_{\D'}(w)\,  |dw|$ on $\D'$
with constant curvature $-4$, that is
\begin{equation*}
\lambda_{\D' }(w) \, |dw| = \frac{1}{2} \, \frac{1}{|w| \, \log
  \frac{1}{|w|}}\, |dw|\, .
\end{equation*} 
Then Lemma \ref{lem:new_sk} and Lemma \ref{lem:hebsing} imply that
\begin{equation*}
\mu(z)\, |dz|=\frac{1}{2} \, \frac{|F'(z)|}{|F(z)| \, \log
  \frac{1}{|F(z)|}}\, |B(z)|\, |dz|= \frac{1}{2} \, \frac{|F'(z)|}{|O(z)| \, \log
  \frac{1}{|F(z)|}}\, |dz|
\end{equation*}
is an SK--metric on $\D$. We
further note that 
\begin{equation*}
\limsup_{z \to z_j} \frac{\mu(z)}{|z-z_j|^{m_j}}<+ \infty
\end{equation*}
 for all $j$. Hence 
$\mu(z) \le \lambda_{max}(z)$ for $z \in \D$ and therefore
\begin{equation*}
 \frac{1}{2} \, \frac{1}{|O(z)| \, \log
  \frac{1}{|F(z)|}} \le \frac{1}{1-|F(z)|^2}\,\quad
 \text{for } z\in \D\, .  \end{equation*}
Letting $z \stackrel{\angle}{\to} \zeta$ yields
\begin{equation*}
\frac{1}{2} \, \frac{1}{|\beta| \, \log
  \frac{1}{|\beta|}} \le \frac{1}{1-|\beta|^2}\,, 
\end{equation*}
violating $\lambda_{\D'} (z) > \lambda_{\D}(z)$ for all $z \in
\D$. Hence $ \angle \lim_{z \to \zeta} |O(z)|=1$ for a.\!\;e.~$\zeta \in \partial
\D$ and therefore $O\equiv c$ for some constant $c$, $|c|=1$. Thus $F$ is a
Blaschke product, as required.

\medskip

It remains to show that $F$ is an indestructible Blaschke product.
We have seen that {\it each} developing map $F$ of $\lambda_{max}(z)\, |dz|$ 
is a Blaschke product. Thus the result follows by  Liouville's theorem
(Theorem \ref{thm:liouville}).
 \hfill{$\blacksquare$}

\medskip

We next consider maximal functions whose critical sets form 
finite and Blaschke sequences respectively.

\begin{theorem}\label{thm:finiteC}
\begin{itemize}
\item[(a)]
Let ${\cal C}$ be a finite sequence of $n$ points in
$\D$.
Then every maximal function  for $\mathcal{C}$ is a finite Blaschke
  product of degree $m=n+1$.
\item[(b)]
If ${\cal C}$
is a Blaschke sequence  in $\D$ then
every maximal function
for ${\cal C}$ is an indestructible Blaschke product which has a finite angular derivative at almost
every point of $\partial \D$.
\end{itemize}
\end{theorem}

{\bf Proof.}
Suppose ${\cal C}$ is a finite or Blaschke sequence. Let $F: \D \to \D$
 be a maximal function for ${\cal C}$ and $\lambda_{max}(z)\,
|dz|= F^*\lambda_{\D}(z)\, |dz|$ be the maximal conformal metric with constant
curvature $-4$ and zero set ${\cal C}$. By
Theorem \ref{thm:0} the maximal function $F$ is an indestructible Blaschke product. 
Now, let $B$ be a Blaschke product with zero set ${\cal C}$.
Then using Lemma \ref{lem:new_sk} we see that
\begin{equation*}
 \lambda(z)\, |dz| :=|B(z)| \, \lambda_{\D}(z)  \, |dz|
\end{equation*}
is a conformal pseudometric on $\D$ with curvature bounded above by $-4$ and zero set
${\cal C}$.
Thus 
\begin{equation*}
\lambda(z) \le \lambda_{max}(z) \quad \text{ for } z \in \D 
\end{equation*}
and consequently
\begin{equation}\label{eq:ab1}
 |B(z)| \le \frac{\lambda_{max}(z)}{\lambda_{\D}(z)}  \quad \text{ for all } z \in \D  \, . 
\end{equation}

\medskip

(a)
If $B$ is a finite Blaschke product,  we deduce from (\ref{eq:ab1}) and the
Schwarz--Pick lemma that
\begin{equation*}
\lim_{z \to \zeta}
\frac{\lambda_{max}(z)}{\lambda_{\D}(z)}=\lim_{z \to \zeta} \frac{|F'(z)|}{1-|F(z)|^2}\, (1-|z|^2) =1 \quad \text{ for
  all } \zeta \in \partial \D \, .
\end{equation*}
Applying Theorem \ref{thm:boundary} (see Subsection \ref{sec:boundary}) shows that $F$ is a finite
Blaschke product.
The branching order of $F$ is clearly $2\!\;n$. Thus,
according to the Riemann--Hurwitz formula, see \cite[p.~140]{For1999}, the
Blaschke product $F$ has degree $m=n+1$.

\medskip

(b)
Suppose now that ${\cal C}$ is a Blaschke sequence. 
Since $\angle \lim_{z \to \zeta} |B(z)|=1$ for a.\!\;e.~$\zeta
\in \partial \D\, , $ it follows  from (\ref{eq:ab1}) and the Schwarz--Pick lemma that

\begin{equation*}
\angle
\lim_{z \to \zeta}
\frac{\lambda_{max}(z)}{\lambda_{\D}(z)}=\angle\lim_{z \to \zeta}
\frac{|F'(z)|}{1-|F(z)|^2}\, (1-|z|^2) =1 \quad \text{ for a.\!\;e. } \zeta \in \partial \D \, .
\end{equation*}
Hence, by Corollary \ref{cor:5}, we conclude that $F$ has a finite angular
derivative at almost every boundary point of $ \D$.
\hfill{$\blacksquare$}

\medskip

As already remarked, Theorem  \ref{thm:finiteC} (a) combined with
Nehari's uniqueness result \cite[Corollary to Theorem 1]{Neh1946} characterizes maximal functions with finite
critical sets. Part (b) of 
Theorem \ref{thm:finiteC} raises the question of whether  Nehari's result
extends to indestructible Blaschke products 
with Blaschke sequences as critical sets and
 a finite angular derivative at almost every boundary
point of $\D$.

\medskip

Combining Theorem \ref{thm:0} and Theorem \ref{thm:finiteC} leads to
Corollary \ref{cor:3} which we are going to prove now.

\medskip

{\bf Proof of Corollary \ref{cor:3}.}
The uniqueness assertion is a direct consequence of
Lemma \ref{lem:gen_max_3} and Liouville's theorem respectively. 

\smallskip
 
For the existence part we consider $\lambda(z)\, |dz|:= f^*\lambda_{\D}(z)\,
|dz|$. Then $\lambda(z)\, |dz|$ is a con\-formal pseudometric on $\D$ with constant
curvature $-4$ and zero set ${\cal C}$. Thus if ${\cal C}^* \subseteq
{\cal C}$, then by Theorem \ref{thm:perron2} there exists   a maximal conformal pseudometric
$\lambda_{max}(z)\, |dz|$ on $\D$ with constant
curvature $-4$ and zero set ${\cal C}^*$. In particular,
\begin{equation*}
\frac{|f'(z)|}{1-|f(z)|^2} = \lambda(z) \le \lambda_{max}(z)=
\frac{|F'(z)|}{1-|F(z)|^2}\, ,\quad z \in \D \,.
\end{equation*}
Now the critical set of the maximal function $F$ is ${\cal C}^*$ and 
$F$ is an indestructible Blaschke product by Theorem \ref{thm:0}. 
If ${\cal C}^*$ is finite then   Theorem
\ref{thm:finiteC} (a) shows that $F$
is a finite Blaschke product.
\hfill{$\blacksquare$}

\bigskip

We conclude this section by giving an intrinsic characterization of maximal conformal pseudo\-metrics
with constant curvature  $-4$ whose zero set form finite and Blaschke sequences. It
might be of interest to see whether this condition
generalizes to all maximal conformal pseudometrics with constant curvature $-4$.

\begin{theorem}
Let ${\cal C}$ be a finite or Blaschke sequence  in $\D$. A conformal pseudometric
$\lambda(z)\, |dz|$ on $\D$ with constant curvature $-4$ and zero set ${\cal
  C}$ is the maximal conformal pseudometric $\lambda_{max}(z)\, |dz|$ on $\D$
with constant curvature $-4$ and zero set ${\cal C}$ if and only if
\begin{equation}\label{eq:boundary2}
\lim_{ r\to 1} \,  \int \limits_0^{2\pi} \log
\frac{\lambda(re^{it})}{\lambda_{\D}(re^{it})}\,dt=0\, .
\end{equation}
\end{theorem}

{\bf Proof.}
Let $B$ be a Blaschke product with zero set ${\cal C}$. Applying
Lemma \ref{lem:new_sk}, we see that \begin{equation*}
 \lambda(z)\, |dz| :=|B(z)| \, \lambda_{\D}(z)  \, |dz|
\end{equation*}
defines a conformal pseudometric on $\D$ with curvature bounded above by $-4$ and zero set
${\cal C}$. Thus
 $\lambda(z) \le \lambda_{max}(z)$ for $z \in \D$ and therefore
\begin{equation*}
 |B(z)| \le \frac{\lambda_{max}(z)}{\lambda_{\D}(z)}  \quad \text{ for all }  z \in \D  \, . 
\end{equation*}
Since $\lim_{ r\to 1} \,  \int \limits_0^{2\pi} \log|B(re^{it})|\, dt=0$, see
\cite[Chapter II, Theorem 2.4]{Gar2007}, 
we deduce that
\begin{equation*}
\lim_{ r\to 1} \,  \int \limits_0^{2\pi} \log
\frac{\lambda_{max}(re^{it})}{\lambda_{\D}(re^{it})}\,dt=0\, ,
\end{equation*}
as desired. 

\smallskip

Conversely, let
$\lambda(z)\, |dz|$ be a conformal pseudometric on $\D$ with constant
curvature $-4$ and zero set ${\cal C}$ which satisfies (\ref{eq:boundary2}).
We now consider the nonnegative subharmonic
function
$z \mapsto \log (\lambda_{max}(z)/\lambda(z))$ on $\D$ and note that
\begin{equation*}
\lim_{ r\to 1} \,  \int \limits_0^{2\pi} \log
\frac{\lambda_{max}(re^{it})}{\lambda(re^{it})}\,dt=  \lim_{ r\to 1} \,  \int \limits_0^{2\pi} \log
\frac{\lambda_{max}(re^{it})}{\lambda_{\D}(re^{it})}\,dt - \lim_{ r\to 1} \,  \int \limits_0^{2\pi} \log
\frac{\lambda(re^{it})}{\lambda_{\D}(re^{it})}\,dt=  0\, .
\end{equation*}
Hence $\lambda(z)\, |dz|=\lambda_{max}(z)\, |dz|$.{\hfill{$\blacksquare$}}

\section{Proof of Theorem \ref{thm:main}}
\label{sec:5}

We now  combine our previous results to complete the proof of Theorem \ref{thm:main}.

\medskip

{\bf Proof of Theorem \ref{thm:main}.}\\
 (a) $\Rightarrow$ (b):  This is Corollary \ref{cor:1}.

\smallskip

(b) $\Rightarrow$ (c): By Remark \ref{rem:little_paley} we have $B' \in {\cal
  A}_1^2$.

\smallskip

(c) $\Rightarrow$ (d):  Again using Remark \ref{rem:little_paley}, we see that
if $\varphi \in {\cal A}_1^2$, then
\begin{equation*}
z \mapsto \int \limits_0^z \varphi(\xi) \, d\xi\, , \quad z \in \D, 
\end{equation*}
belongs to $H^2$ and therefore to the
Nevanlinna class ${\cal N}$.

\smallskip

\label{page5}
(d) $\Rightarrow$ (a):  Let $\varphi \in {\cal N}$. Then $\varphi=
\varphi_1/\varphi_2$ is the quotient of two 
analytic self--maps of $\D$, see for instance
 \cite[Theorem 2.1]{Dur2000}. W.\!\;l.\!\;o.\!\;g.~we may assume $\varphi_2$ is zerofree. Differentiation of $\varphi$ yields
\begin{equation*}
\varphi'(z)=\frac{1}{\varphi_2(z)^2} \, \big(
 \varphi_1'(z) \, \varphi_2 (z)- \varphi_1(z)\, \varphi_2'(z)
\big)\, .
\end{equation*}
Since $\varphi'_1, \varphi_2' \in {\cal A}_1^2$ and ${\cal A}_1^2$ is a vector
space,   it follows that the function $\varphi_1' \, \varphi_2 - \varphi_1\, \varphi_2'$
belongs to ${\cal A}_1^2$. Thus Theorem \ref{thm:d} ensures the existence of a
solution $u : \D \to \R$ to 
\begin{equation*}
\Delta u= \left| \varphi'(z) \right|^2 \, e^{2u}\, .
\end{equation*}
By Liouville's theorem 
\begin{equation*}
u(z) =\log \left(\frac{1}{|\varphi'(z)|}\, \frac{|f'(z)|}{1-|f(z)|^2}   \right)\,, 
\quad z \in \D\, ,
\end{equation*}
for some analytic self--map  $f$ of $\D$. Hence the critical set of $\varphi$
coincides with the critical set of $f$, as required.\hfill{$\blacksquare$}

\medskip

We wish to point out that Theorem \ref{thm:main} 
characterizes 
the class of analytic functions whose critical sets
agree with the critical sets of the class of bounded analytic functions.
In fact, the critical set ${\cal
  C}_{g}$ of an analytic function $g: \D \to \C$ coincides with the critical
set ${\cal C}_f$ of an analytic function $f: \D \to \D$ if and only if
$g'=\varphi_1\, \varphi_2$ for some function $\varphi_1 \in {\cal A}_1^2$ and a
nonvanishing analytic function $\varphi_2: \D \to \C$.

\medskip

A final remark is that  the list of equivalent statements in Theorem \ref{thm:main} can be
extended.

\begin{remark}\label{rem:(e)}
Let $ (z_j)$ be a sequence in $\D$. Then the following statements are
equivalent.
\begin{itemize}
 \item[(a)] There is an analytic self--map of $ \D$ with critical set
   $(z_j)$.
\item[(e)] There is a meromorphic function $g$ on $\D$ with critical set
  $(z_j)$, where $g$ is the quotient of two
  analytic self--maps of $\D$.
\end{itemize}
\end{remark}

\smallskip

{\bf Proof.}
We adapt the idea already used in the proof of implication (d) $\Rightarrow$ (a) of
Theorem \ref{thm:main}.
Let $(z_j)$ be the critical set of the meromorphic function $g= \varphi_1/\varphi_2$, where $\varphi_1$ and $\varphi_2$ are analytic
self--maps of $\D$. We  may  assume that $\varphi_1$ and  $\varphi_2$ have
no common zeros. Since $g'(z) \, \varphi_2(z)^2$ belongs to ${\cal A}_1^2$, 
Theorem \ref{thm:d} guarantees a solution $u: \D \to \R$ to
\begin{equation*}
\Delta u= \left| g'(z) \, \varphi_2(z)^2 \right|^2 \, e^{2u}\, .
\end{equation*}
 Now Liouville's theorem shows  that 
\begin{equation*}
u(z) =\log \left(\frac{1}{|g'(z)\, \varphi_2(z)^2|}\, \frac{|f'(z)|}{1-|f(z)|^2}   \right)\,, 
\quad z \in \D\, ,
\end{equation*}
for some analytic self--map  $f$ of $\D$. Thus the critical set of $f$ coincides
with the critical set of $g$ and the zero set of $\varphi_2^2$. 
Now Lemma \ref{lem:subset}  below tells us that every subset
of a critical set of a bounded analytic function is again the critical set of
another bounded analytic function. Hence there exists an analytic
self--map of $\D$ with critical set $(z_j)$, as desired.
\hfill{$\blacksquare$}

\medskip

\begin{lemma}\label{lem:subset}
Every subset of a critical set of a bounded analytic function on $\D$ is the
critical set of another  bounded analytic function on $\D$.
\end{lemma}

\medskip

{\bf Proof.}
It suffices to consider nonconstant analytic self--maps of $\D$.

\smallskip

The statement, using Theorem \ref{thm:main}, follows from the corresponding
result about the zero sets
of functions in ${\cal A}_1^2$, see \cite[Theorem 7.9]{Hor1974} and
\cite[Corollary 4.36]{HKZ}.

\smallskip

Here, we would like to give an alternative proof.
Suppose $f$ is a nonconstant
analytic self--map  of  $\D$ with critical set ${\cal C}$. Now let $ {\cal
  C}^*$ be a subsequence of ${\cal C}$.
By Theorem \ref{thm:perron2} there exists a regular conformal
pseudometric on $\D$ with constant curvature $-4$ and  zero set ${\cal C}^*$.
Applying Liouville's theorem gives the desired analytic self--map $g$ of $\D$ with critical
set ${\cal C}^*$.~\hfill{$\blacksquare$}

\hspace{16mm}

\hspace{0.9cm}\begin{minipage}{8cm}
Daniela Kraus\\
Department of Mathematics\\
University of W\"urzburg\\
97074 W\"urzburg\\
Germany\\
dakraus@mathematik.uni-wuerzburg.de

\end{minipage}

\begin{thebibliography}{99}
\addcontentsline{toc}{section}{References}

\fontsize{11pt}{0.pt}
\selectfont
{
\bibitem{Ahl1938} L.~Ahlfors, An extension of Schwarz's lemma, {\it
    Trans.~Am.~Math.~Soc.}~{\bf 43} (1938), 359--364.

\bibitem{Aub1998} T.~Aubin, {\it Some nonlinear problems in Riemannian geometry},
  Springer 1998.

\bibitem{Avi1986} P.~Aviles, Conformal complete metrics with prescribed
  non--negative Gaussian curvature in $\R^2$, {\it Invent.~Math.}~{\bf 83}
  (1986), 519--544.


\bibitem{BM2007} A.~Beardon and D.~Minda, The Hyperbolic Metric and Geometric
  Function Theory, in: S.~Ponnusamy, T.~Sugawa and M.~Vuorinen (editors), {\it
    Quasiconformal Mappings and their Applications}, Narosa 2007.

\bibitem{BK1986} J.~Bland and M.~Kalka, Complete metrics conformal to the
  hyperbolic disc, {\it Proc. Am.~Math.~Soc.} {\bf 97} (1986), no.~1, 128--132. 






\bibitem{Bla1915} W.~Blaschke, Eine Erweiterung des Satzes von Vitali \"uber
  Folgen analytischer Funktionen, {\it S.--B.~S\"achs.~Akad.~Wiss.~Leipzig
    Math.--Natur.\,Kl.} {\bf 67} (1915), 194--200.


\bibitem{Bie16} L.~Bieberbach,  $\Delta u=e^u$ und die
automorphen Funktionen, {\it Math.~Ann.}~{\bf 77} (1916), 173--212.

\bibitem{Chang2004} S.~Y.~A.~Chang,
{\it  Non--linear elliptic equations in conformal geometry},
European Mathematical Society 2004.

\bibitem{Chen2001}
H.~Chen,  On the Bloch constant, in:
 N.~Arakelian et al. (editors),  Approximation, complex analysis, and potential
 theory,  Kluwer Academic Publishers (2001),
 129--161.



\bibitem{CN1991} K.~S.~Cheng and W.~M.~Ni, On the structure of the conformal
  Gaussian curvature equation on $\R^2$, {\it Duke
    Math.\!\;J.}~{\bf 62 }  (1991), no.~3, 721--737.


\bibitem{CW94} K.~S.~Chou and T.~Wan, Asymptotic radial symmetry
for solutions of $\Delta u+e^u=0$ in a punctured disc, {\it Pacific
  J.~Math.}~{\bf 163} (1994), no.~2, 269--276. 

\bibitem{CW95} K.~S.~Chou and T.~Wan,  Correction to  ``Asymptotic radial symmetry
for solutions of $\Delta u+e^u=0$ in a punctured disc'', 
{\it Pacific J.~Math.}~{\bf 171} (1995), no.~2, 589--590.
 
\bibitem{Col1985} P.~Colwell, {\it Blaschke products}, The University of Michigan
  Press 1985.

\bibitem{Cou1968} R.~Courant and D.~Hilbert, {\it Methoden der Mathematischen Physik II}, Springer 1968.


\bibitem{Dur1969}
P.~Duren,
On the Bloch-Nevanlinna conjecture, {\it Colloq.~Math.}~{\bf 20} (1969), 295--297.


\bibitem{Dur2000} P.~Duren, {\it Theory of $H^p$ spcaes}, Dover 2000.


\bibitem{For1999} O.~Forster, {\it Lectures on Riemann surfaces}, Springer 1999.


\bibitem{Gar2007} J.~B.~Garnett, {\it Bounded analytic functions}, revised
  first edition, Springer 2007.


\bibitem{GT} D.~Gilbarg and  N.~S.~Trudinger, {\it Elliptic partial differential
equations of second order}, Springer 2001.





\bibitem{HKZ} H.~Hedenmalm, B.~Korenblum and K.~Zhu, {\it Theory of Bergman
    Spaces}, Springer 2000. 


\bibitem{Hei62} M.~Heins, 
On a class of conformal metrics, {\it Nagoya Math.\!\;J.} {\bf 21} (1962),
    1--60.

\bibitem{Hei86} M.~Heins,
Some characterizations of finite Blaschke products of positive degree, {\it
  J.\!\;Anal. Math.}~{\bf 46} (1986), 162--166.


\bibitem{Hor1974} C.~Horowitz, Zeros of functions in the Bergman spaces, 
{\it Duke Math.\!\;J.}~{\bf 41} (1974), 693--710.

\bibitem{Hor1977} C.~Horowitz, Factorization theorems for functions in the
  Bergman spaces, {\it Duke Math.\!\;J.}~{\bf 44} (1977), 201--213.



\bibitem{HT92} D.~Hulin and M.~Troyanov,  Prescribing curvature on open
  surfaces, {\it Math.~Ann.}~{\bf 293} (1992), no.~2, 277--315.




\bibitem{Jen1899} J.~Jensen,  Sur un nouvel et important th\'{e}or\`{e}me de la th\'{e}orie des fonctions,
{\it Acta Math.}~{\bf 22} (1899), 359--364.



\bibitem{Kaz1985} J.~Kazdan, Prescribing the curvature of a Riemannian
  manifold, {\it CMBS Regional Conf.~Ser.~in Math.}~{\bf 57} (1985). 

\bibitem{KY1993} M.~Kalka and D.~Yang, On conformal deformation of nonpositive
  curvature on noncompact surfaces, {\it Duke
    Math.\!\;J.}~{\bf 72} (1993), no.~2, 405--430.


\bibitem{KL2007} L.~Keen and N.~Lakic, {\it Hyperbolic Geometry from a Local
    Viewpoint}, Camb. Univ. Press 2007. 

\bibitem{Kor1975} B.~Korenblum, An extension of the Nevanlinna theory, {\it Acta
  Math.}~{\bf 135} (1975), 187--219. 


\bibitem{Kra2011b} D.~Kraus,  Zero sets of weigthed Bergman
    spaces and the Gauss curvature equation, preprint.

\bibitem{KR} D.~Kraus and O.~Roth,  Critical points of inner functions, nonlinear partial differential equations, and an extension of Liouville's theorem,
{\it J.~London Math.~Soc}.  {\bf 77} (2008), no.~1, 183--202.

\bibitem{KR2008}  D.~Kraus and O.~Roth, 
Conformal Metrics, in:~\textit{Topics in Modern Function Theory}, Ramanujan
Math.~Soc., 41 pp., to appear.


\bibitem{KRR06} D.~Kraus, O.~Roth and St.~Ruscheweyh,
A boundary version of Ahlfors' Lemma, locally complete  conformal metrics
and conformally invariant reflection  principles for analytic maps,
{\it J.~Anal.~Math.}~{\bf 101} (2007), 219--256.


 \bibitem{BF2008} I.~Laine, Complex differential equations, in: F.~Battelli
   and M.~Fe\v{c}kan (editors),  {\it Handbook of
     differential equations:~Ordinary differential equations, Vol.~IV},
   Elsevier 2008.

\bibitem{Lio1853} J.~Liouville,  
Sur l'\'equation aux diff\'erences
  partielles $\frac{d^2 \log \lambda}{du dv}\pm \frac{\lambda}{2
  a^2}=0$,  {\it J.\!\;de Math.}~{\bf 16} (1853), 71--72.

\bibitem{MT2002} R.~Mazzeo and M.~Taylor, Curvature and Uniformization, {\it
    Isr.~J.~Math.}~{\bf 130} (2002), 323--346.



\bibitem{Min1987} D.~Minda, The strong form of Ahlfors' lemma, {\it Rocky
    Mt.~J.~Math.}~{\bf 17} (1987),  no.~3, 457--461.


\bibitem{Min} D.~Minda, Conformal metrics, unpublished notes.

\bibitem{Mos1973} J.~Moser,
On a nonlinear problem in differential geometry,
{\it Dynamical Syst., Proc. Sympos. Univ. Bahia, Salvador 1971} (1973), 273--280.


\bibitem{Neh1946} Z.~Nehari, A generalization of Schwarz' lemma, {\it Duke
    Math.\!\;J.}~{\bf 14} (1947), 1035--1049.


\bibitem{Nevs1922} F.~and R.~Nevanlinna, \"Uber die Eigenschaften analytischer
  Funktionen in der Umgebung einer singul\"aren Stelle oder Linie, {\it Acta
    Soc.~Sci.~Fenn.}~{\bf 50} (1922), no.~5, 46 pp.


\bibitem{Ni1989} W.--M.~Ni, Recent progress on the elliptic equation $\Delta
  u + K\, e^{2u}=0$ on $\R^2$, {\it Rend. Semin. Mat., Torino Fasc. Spec.}~1989, 1--10.   

\bibitem{Nit57} J.~Nitsche,  \"Uber die isolierten
Singularit\"aten der L\"osungen von $\Delta u=e^{u}$,
{\it Math.~Z.}~ {\bf 68} (1957), 316--324. 

\bibitem{Per1923} O.~Perron, Eine neue Behandlung der Randwertaufgabe f\"ur
$\Delta u=0$, {\it Math.~Z.}~{\bf 18} (1923), 42--54.


\bibitem{Pri1956} I.~L.~Priwalow, {\it Randeigenschaften analytischer
Funktionen}, Deutscher Verlag der Wissenschaften 1956.



\bibitem{Roy1986} H.~L.~Royden, The Ahlfors--Schwarz lemma: the
  case of equality, {\it J.~Anal.~Math.}~{\bf 46} (1986), 261--270.



\bibitem{Seip1994} K.~Seip, On a theorem of Korenblum, {\it Ark.~Math.}~{\bf 32} (1994), 237--243.

\bibitem{Seip1995} K.~Seip, On Korenblum's density condition for zero
  sequences of $A^{- \alpha}$, {\it J.~Anal.~Math.}~ {\bf 67} (1995), 307--322.



\bibitem{Sha1993} J.~H.~Shapiro, {\it Composition operators and classical
    function theory}, Springer 1993.


\bibitem{Smi1986} S.~J.~Smith, {\it On the uniformization of the $n$--punctured
    disc}, Ph.\!\;D,~Thesis, University of New England, 1986.


\bibitem{Ste2005} K.~Stephenson, {\it Introduction to circle packing: the theory of discrete analytic
functions}, Camb.~Univ.~Press 2005.


\bibitem{Str2005} M.~Struwe,
A flow approach to Nirenberg's problem, {\it
Duke Math.\!\;J.}~{\bf 128} (2005), no.\!\;1, 19--64.



\bibitem {Wal1950} J.~Walsh, {\it The location of critical points of analytic
    and harmonic functions}, Am.~Math.~Soc.~1950.

\bibitem{WP79} Q.~Wang and J.~Peng, On critical points of finite Blaschke
  products and the equation $\Delta u=e^{2u}$, {\it Kexue Tongbao}~{\bf 24} (1979),
   583--586 (Chinese).

\bibitem{Yam1988} A.~Yamada, Bounded analytic functions and metrics of
  constant curvature on Riemann surfaces, {\it Kodai Math.\!\;J.}~{\bf 11}
  (1988), no.~3, 317--324.


\bibitem{Yau1975}  S.~T.~Yau, Harmonic functions on complete Riemannian
  manifolds, {\it Comm.~Pure Appl.~Math.}~{\bf 28} (1975), 201--228.

\bibitem{Yau1978} S.~T.~Yau, A general Schwarz lemma for K\"ahler manifolds, 
  {\it Am.~J.~Math.}~{\bf 100} (1978), 197--203. 

\bibitem{Z96}
S.~Zakeri, On critical points of proper holomorphic maps on the unit disk,
{\it Bull.~London Math.~Soc.}~{\bf 30} (1996), no.~1, 62--66.
}


\end{thebibliography}
\end{document}